\documentclass[11pt]{amsart}
\usepackage{amsmath,amsthm,color,amssymb, amscd, amsfonts,url,multicol}
\usepackage[notcite,notref,final]{showkeys}
\usepackage{pdflscape}

\usepackage[dvipsnames]{xcolor}

\usepackage{tikz}

\usetikzlibrary{calc} 

\usepackage{amsxtra}
\usepackage{amssymb}
\usepackage{mathtools}
\usepackage{mathrsfs}  
\theoremstyle{plain}
\newcommand {\sectionnew}[1]{\section{#1}}

\newtheorem{theorem}{Theorem}[section]
\newtheorem{lemma}[theorem]{Lemma}
\newtheorem{sublemma}[theorem]{Sublemma}
\newtheorem{definition-theorem}[theorem]{Definition-Theorem}
\newtheorem{definition-lemma}[theorem]{Definition-Lemma}
\newtheorem{proposition}[theorem]{Proposition}
\newtheorem{corollary}[theorem]{Corollary}
\newtheorem{conjecture}[theorem]{Conjecture}
 \theoremstyle{definition}
\newtheorem{definition}[theorem]{Definition}
\newtheorem{example}[theorem]{Example}
\newtheorem{remark}[theorem]{Remark}
\newtheorem{notation}[theorem]{Notation}
\newtheorem{routine}[theorem]{Routine}
\newtheorem{problem}[theorem]{Problem}
\newtheorem{hypothesis}[theorem]{Hypothesis}
\newtheorem{question}[theorem]{Question}
\newcommand \bth[1] { \begin{theorem}\label{t#1} }
\newcommand \ble[1] { \begin{lemma}\label{l#1} }
\newcommand \bsubl[1] { \begin{sublemma}\label{sl#1} }

\newcommand \bdl[1] { \begin{definition-lemma}\label{dl#1} }
\newcommand \bpr[1] { \begin{proposition}\label{p#1} }
\newcommand \bco[1] { \begin{corollary}\label{c#1} }
\newcommand \bde[1] { \begin{definition}\label{d#1}\rm }
\newcommand \bex[1] { \begin{example}\label{e#1}\rm }
\newcommand \bre[1] { \begin{remark}\label{r#1}\rm }
\newcommand \bcj[1] { \begin{conjecture}\label{j#1}\rm }
\newcommand \bqu[1]  { \begin{question}\label{n#1}\rm }
\newcommand \bnota[1] { \begin{notation}\label{n#1}\rm }
\newcommand \bro[1] { \begin{routine}\label{n#1}\rm }
\newcommand \bpb[1] { \begin{problem}\label{n#1}\rm }
\newcommand \bhy[1] { \begin{hypothesis}\label{n#1}\rm }
\renewcommand {\eth} { \end{theorem} }
\newcommand {\ele} { \end{lemma} }
\newcommand {\esubl} { \end{sublemma} }

\newcommand {\edl}{ \end{definition-lemma} }
\newcommand {\epr} { \end{proposition} }
\newcommand {\eco} { \end{corollary} }
\newcommand {\ede} { \end{definition} }
\newcommand {\eex} { \end{example} }
\newcommand {\ere} { \end{remark} }
\newcommand {\ecj} { \end{conjecture} }
\newcommand {\equ} {\end{question}}
\newcommand {\enota} { \end{notation} }
\newcommand {\ero} { \end{routine} }
\newcommand {\epb} { \end{problem} }
\newcommand {\ehy} { \end{hypothesis} }
\newcommand \thref[1]{Theorem \ref{t#1}}
\newcommand \dlref[1]{Definition-Lemma \ref{dl#1}}
\newcommand \leref[1]{Lemma \ref{l#1}}
\newcommand \prref[1]{Proposition \ref{p#1}}

\newcommand \deref[1]{Definition \ref{d#1}}

\newcommand \lb[1]{\label{#1}}

\makeatletter              % This sequence of commands will
\let\c@equation\c@theorem  % incorporate equation numbering
\makeatother
\numberwithin{equation}{section}

\def \Cset {{\mathbb C}}

\def \Zset {{\mathbb Z}}
\def \Nset {{\mathbb N}}
\def \Qset {{\mathbb Q}}

\def \AA  {{\mathcal{A}}}           %mathcal

\def \CC {{\mathcal{C}}}

\def \PP {{\mathcal{P}}}

\def \hb {{\hbar}}
\def \al {\alpha}
\def \be {\beta}

\def \ga {\gamma}

\def \mt  {\mapsto}
\def \hra {\hookrightarrow}

\def \rcor {\rangle}
\def \lcor {\langle}
\def \del {\partial}

\def \ol {\overline}
\def \wt {\widetilde}
\def \wh {\widehat}

\def \id { {\mathrm{id}} }

\def \rank { {\mathrm{rank}} }

\def \nn  {\mathfrak{n}}
\def \mm  {\mathfrak{m}}

\newcommand{\kk}{\Bbbk}

\DeclareMathOperator \Aut { {\mathrm{Aut}} }
\DeclareMathOperator \Der { {\mathrm{Der}} }

\DeclareMathOperator \maxSpec { {\mathrm{maxSpec}}}

\DeclareMathOperator \ad { {\mathrm{ad}} }

\DeclareMathOperator \tr { {\mathrm{tr}} }

\DeclareMathOperator \Irr { {\mathrm{Irr}} }

\renewcommand \Im { {\mathrm{Im}} }

\newcommand \red {\textcolor{black}}
\newcommand \blue {\textcolor{black}}
\newcommand \green {\textcolor{black}}

\newcommand{\pcoor}[1]{%
  \begingroup\lccode`~=`: \lowercase{\endgroup
  \edef~}{\mathbin{\mathchar\the\mathcode`:}\nobreak}%
  [% opening symbol
  \begingroup
  \mathcode`:=\string"8000
  #1%
  \endgroup
  ]% closing symbol
}

\input xy
\xyoption{all}

\begin{document}
%%%%%%%%%%%%%%%%%%%%%%%%%%%%%%%%%%%%%%%%%%%%%%%%%%%%%%%%%%%%%%%%%%%%%%%%%%%
%%%%%%%%%%%%%%%%%%%%%%    Title    %%%%%%%%%%%%%%%%%%%%%%%%%%%%%%%%%%%%%%%%

\title[Poisson geometry of PI  3-dimensional Sklyanin algebras]
{Poisson geometry of PI 3-dimensional \\ Sklyanin algebras}
\author[Chelsea Walton]{Chelsea Walton}
\address{
\green{Department of Mathematics, 
The University of Illinois Urbana-Champaign,
Urbana, IL 61801,~USA}
}
\email{\green{notlaw@illinois.edu}}
\author[Xingting Wang]{Xingting Wang}
\address{
\green{Department of Mathematics \\
Howard University \\
Washington, DC 20059, 
USA}
}
\email{\green{xingting.wang@howard.edu}}
\author[Milen Yakimov]{Milen Yakimov}
\address{
Department of Mathematics \\
Louisiana State University \\
Baton Rouge, LA 70803,
USA
}
\email{yakimov@math.lsu.edu}
\thanks{\blue{Walton was partially supported by NSF grant \#1550306 and a research fellowship from the Alfred P. Sloan Foundation.
Yakimov was supported by NSF grant \#1601862 and Bulgarian Science Fund grant H02/15.
Wang was partially supported by an AMS-Simons travel grant.}}
%\date{}

\begin{abstract} 
We give the 3-dimensional Sklyanin algebras $S$ that are module-finite over their center $Z$ the structure of a
Poisson $Z$-order (in the sense of Brown-Gordon). We show that the induced Poisson bracket on $Z$ is non-vanishing and is induced by an explicit potential. 
The $\Zset_3 \times \kk^\times$-orbits 
of symplectic cores of the Poisson structure are determined (where the group acts on $S$ by algebra automorphisms). In turn, 
this is used to analyze the finite-dimensional quotients of $S$ by central annihilators: there are 3 distinct isomorphism classes of such
quotients in the case $(n,3) \neq 1$ and 2 in the case $(n,3)=1$, where $n$ is \green{the} order of the elliptic curve automorphism associated to $S$. 
The Azumaya locus of $S$ is determined, extending results of Walton for~the~case~$(n,3)=1$.
\end{abstract}

\maketitle

%\bibliographystyle{abbrv}  

%\setcounter{tocdepth}{2} \tableofcontents

%%%%%%%%%%%%%%%%%%%%   Introduction   %%%%%%%%%%%%%%%%%%%%%%%%%%%%%%%%%%%%%%%%
\sectionnew{Introduction}
\lb{Intro}
Throughout the paper, $\kk$ will denote an algebraically closed field of characteristic $0$. In 2003, Kenneth Brown and Iain Gordon \cite{BrownGordon} introduced the notion of a {\it Poisson order} in order to provide a framework for studying the representation theory of algebras that are module-finite over their center, with the aid of Poisson geometry. A {\it Poisson $C$-order} is a finitely generated $\kk$-algebra $A$ that is module-finite over a central subalgebra $C$ so that there is a $\kk$-linear map from $C$ to the space of derivations of $A$ that imposes on $C$ the structure of a Poisson algebra (see Definition~\ref{dPois-order}). Towards studying the irreducible representations of such $A$, Brown and Gordon introduced  a {\it symplectic core stratification} of the affine Poisson variety $\maxSpec C$ which is a coarsening of the symplectic foliation of $C$ if $\kk = \Cset$. Now the surjective map from the set of irreducible representations $I$ of $A$ to their central annihilators $\mathfrak{m} = \green{\text{Ann}_A(I) \cap C}$ have connected fibers along symplectic cores of $C$. In fact, Brown and Gordon proved remarkably that for $\mm$ and $\nn$ in the same symplectic core of $C$ we have an isomorphism of the corresponding finite-dimensional, central quotients of $A$,
\[
A/(\mm A) \cong  A/(\nn A).
\]
Important classes of noncommutative algebras arise as Poisson orders, and thus have representation theoretic properties given by symplectic cores. Such algebras include many quantum groups at roots of unity (as shown in \cite{DKP,DL,DP}) and symplectic reflection algebras \cite{EtingofGinzburg} (as shown in \cite{BrownGordon}). Our goal in this work is to show that this list includes another important class of noncommutative algebras, the {\it 3-dimensional Sklyanin algebras} that are module-finite over their center, and to apply this framework to the structure and representation theory of these algebras. Previously Poisson geometry was not utilized for these types of representation theoretic questions; cf. \cite{deLaetLeBruyn, ReichWalton, Walton}.

\smallbreak Three-dimensional Sklyanin algebras (Definition \ref{dSkly3}) arose in the 1980s through Artin, Schelter, Tate, and van den Bergh's classification of noncommutative graded analogues of commutative polynomial rings in 3 variables \cite{AS, ATV1}; these algebras were the most difficult class to study as they are quite tough to analyze with traditional Gr\"{o}bner basis techniques (see, e.g., \cite{IS}). So, projective geometric data was assigned to the algebras $S$, namely an elliptic curve $E \subset \mathbb{P}^2$, an invertible sheaf $\mathcal{L}$ on $E$, and an automorphism $\sigma$ of $E$, in order to analyze the ring-theoretic and homological behaviors of $S$ \cite{ATV1}. It was shown that there always exists a central regular element $g$ of $S$ that is homogeneous of degree 3. Moreover, $S/gS$ is isomorphic to a {\it twisted homogeneous coordinate ring} $B:=B(E, \mathcal{L}, \sigma)$, a noncommutative version of the homogeneous coordinate ring $\bigoplus_{i \geq 0} H^0(E, \mathcal{L}^{\otimes i})$. 

\smallbreak Along with numerous good properties of $S$  (e.g., being a Noetherian domain of polynomial growth and global dimension 3) that were established in \cite{ATV1} (often by going through~$B$), it was shown in \cite{ATV2} that $S$ and $B$ are module-finite over their centers precisely when $\sigma$ has finite order. Now take $$n:= |\sigma|.$$ In this case, the structure of the center $Z$ of $S$ was determined in \cite{AST, SmithTate}: it is generated by $g$ and algebraically independent variables $z_1, z_2, z_3$ of degree $n$, subject to one relation $F$ of degree $3n$. Moreover in this case, $S$ is module-finite over its center if and only if it is a polynomial identity (PI) algebra of PI degree $n$ (see \cite[Corollary~3.12]{Walton}). The maximum dimension of the irreducible representations of $S$ is $n$ as well by \cite[Proposition~3.1]{BrownGoodearl}.

\smallbreak Recall that the Heisenberg group $H_3$ is the group of upper triangular $3 \times 3$-matrices with entries in $\Zset_3$ and 1's on the diagonal. 
It acts by graded automorphisms on $S$ in such a way that the generating space $S_1$ is the standard 3-dimensional representation of $H_3$. On the other hand, the center of $S$ was described in detail 
by Smith and Tate \cite{SmithTate} in terms of so-called {\it good bases} $\{x_1, x_2, x_3 \}$, see Section~\ref{3DSkly}. There is a good basis $\{x_1, x_2, x_3 \}$
that is cyclically permuted by one of the generators $\tau$ of $H_3$. In our study of the geometry of Poisson orders on $S$, we will make an essential use 
of the action of the group $\Sigma:= \mathbb{Z}_3 \times \Bbbk^\times$ on $S$ by graded algebra automorphisms. Here, $\mathbb{Z}_3 = \langle \tau \rangle$ acts on $x_1,x_2, x_3$ as above and $\Bbbk^\times$ acts by simultaneously rescaling the generators.

\smallbreak Our first theorem is as follows.

\bth{intro}
Let $S$ be a 3-dimensional Sklyanin algebra that is module-finite over its center $Z$ and retain the notation above.  Suppose that $\{z_1, z_2, z_3\}$ are of the form \eqref{good-u}, i.e., given in terms of a good basis of the generating space $S_1 \cong B_1$ (see Definition~\ref{dgood}). Then:

\begin{enumerate}
\item $S$ admits the structure of a $\Sigma$-equivariant Poisson $Z$-order for which the induced Poisson structure on $Z$ has a nonzero Poisson bracket.
\item The formula for the Poisson bracket on $Z$ is determined as follows:
$$\quad \quad \quad  \{z_1,z_2\} = \partial_{z_3} F, \quad \{z_2,z_3\} =  \partial_{z_1} F, \quad \{z_3,z_1\} =  \partial_{z_2} F,$$
with $g$ in the Poisson center of the Poisson algebra $Z$.
\end{enumerate}
\eth

\smallbreak
Take 
$$Y := \maxSpec(Z)  = \mathbb{V}(F) \subset \mathbb{A}^4_{(z_1, z_2, z_3, g)},$$ 
and let $Y^{\textnormal{sing}}$ and $Y^{\textnormal{smooth}} = Y \setminus Y^{\textnormal{sing}}$ be its singular and smooth loci, respectively. 
\thref{intro}(2) turns $Y$ into a singular Poisson variety and the group $\Sigma$ acts on it and on $Z$ by Poisson automorphisms. 
By Lemma \ref{lgoodprelim} (see also Notation~\ref{nnotation-x}),  there exists a good basis of $S_1 ~(\cong B_1)$ so that the subgroup $\mathbb Z_3$ acts on $Z(S)$ 
by cyclically permuting $z_1,z_2,z_3$ and fixing $g$. {\it We work with such a good basis for the rest of the paper.}
The $\kk^\times$-action on 
$Y \subset \mathbb{A}^4_{(z_1, z_2, z_3, g)}$ is given by dilation
\begin{equation}
\label{dilation}
\beta \cdot (z_1, z_2, z_3, g) = (\beta^n z_1, \beta^n z_2, \beta^n z_3, \beta^3 g), \quad \beta \in \kk^\times.
\end{equation}

We prove in Lemma~\ref{lYsing} that
$ Y^{\textnormal{sing}} = \{\underline{0}\}$ if $(n,3) =1$, and is equal to the union of three dilation-invariant curves meeting at $\{\underline{0}\}$ if $(n,3)  \neq 1$ as \green{shown in Figure~1} below. In the second case, the curves have the form $z_i = \al g^{n/3},z_{i+1}=z_{i+2}=0$, $i=1,2,3$ (indices taken modulo 3), for
$\al \in \kk^\times$, that is, each of them is the 
closure of a single dilation orbit. We will denote these curves by $C_1, C_2, C_3$. 
The action of $\Zset_3 \subset \Sigma$ cyclically permutes them. Denote the slices
\[
Y_\gamma := \mathbb{V}(F, g - \gamma) \subset Y \quad \mbox{with} \quad \gamma \in \kk.
\]

Now pertaining to the representation theory of $S$, let $\mathcal{A}$ denote the Azumaya locus of $S$, which is the subset of $Y$ that consists of central annihilators of irreducible representations of maximal dimension; it is open and dense in $Y$ by \cite[Theorem~III.1.7]{BrownGoodearl:book}. Our second theorem is as follows.

\bth{intro2} Retain the notation from above.
\begin{enumerate}
\item The Azumaya locus $\mathcal{A}$ of $S$ is equal to
$Y^{\textnormal{smooth}}.$
\smallskip
\item Each slice $Y_\gamma := \mathbb{V}(F, g - \gamma)$ is a Poisson subvariety of $Y$ with $\gamma \in \kk$, and  
\begin{enumerate}

\item $(Y^{\textnormal{sing}})_\gamma = (Y_\gamma)^{\textnormal{sing}} = Y_\gamma \cap Y^{\textnormal{sing}}$ is the union of the symplectic points of $Y_\gamma$: namely, 
\begin{center}
\smallskip
$\begin{cases} Y_0^{\textnormal{sing}} = \{\underline{0}\} \text{ and } Y_{\gamma \neq 0}^{\textnormal{sing}}= \varnothing, &\text{ if } (n,3) =1;\\
 Y_0^{\textnormal{sing}} = \{\underline{0}\} \text{ and } Y_{\gamma \neq 0}^{\textnormal{sing}} \text{ is the union of 3 distinct points}, &\text{ if }(n,3) \neq 1.
\end{cases}$
\smallskip
\end{center}
\item $Y_{\gamma} \backslash Y_\gamma^{\textnormal{sing}}$ is a symplectic core of $Y$. 
\end{enumerate}

\smallskip

\item The $\Sigma$-orbits of the symplectic cores of $Y$ are 
\begin{enumerate}
\item $Y \backslash Y_0$ if $(n,3)=1$ \quad and \quad  $Y \backslash (Y_0 \cup C_1 \cup C_2 \cup C_3)$ if $(n,3) \neq 1$,
\item $(C_1 \cup C_2 \cup C_3) \backslash \{ \underline{0} \}$ if $(n,3) \neq 1$,
\item $Y_0 \backslash \{ \underline{0} \}$, and
\item $\{ \underline{0} \}$.
\end{enumerate}
These sets define a partition of $Y$ by smooth locally closed subsets with 4 strata in the case $(n,3) \neq 1$ and 3 strata in the case $(n,3)=1$.

\smallskip

\item For each $\mm, \nn \in \maxSpec Z =Y$ that are in the same stratum for the partition in part \textnormal{(3)},
the corresponding central quotients of $S$ are isomorphic: 
\[
S/(\mm S) \cong S/(\nn S).
\]
The central quotients for the strata \textnormal{(3a)} and \textnormal{(3c)} are isomorphic to the matrix algebra $M_n(\kk)$ (Azumaya case). 
In the case \textnormal{(3d)} the central quotient is a local algebra (trivial case). (We refer to case \textnormal{(3b)} as the intermediate case.)
\end{enumerate}
\eth

\vspace{.1in}

\hspace{-.2in}
\begin{tikzpicture}[circle dotted/.style={dash pattern=on .05mm off 1mm,
                                         line cap=round}]
\node (A) at (1.2, 0){} ;
\node (B) at (1.2, 4.5) {};
\node (C) at (9.2, 0) {};
\node (D) at (9.2, 4.5) {};

\draw[<->, dashed]
  (A) edge node[left]{$\scriptstyle g$} (B) (C) edge node[left]{$\scriptstyle g$} (D);

\path[fill, draw, color=white, inner color=black, fill opacity=0.2] (1,0) -- (5,0) -- (6,2) -- (2,2) -- (1,0);
\path[fill, draw, color=white, inner color=black, fill opacity=0.2] (1,2.5) -- (5,2.5) -- (6,4.5) -- (2,4.5) -- (1,2.5);
\path[fill, draw, color=white, inner color=black, fill opacity=0.2] (9,0) -- (13,0) -- (14,2) -- (10,2) -- (9,0);
\path[fill, draw, color=white, inner color=black, fill opacity=0.2] (9,2.5) -- (13,2.5) -- (14,4.5) -- (10,4.5) -- (9,2.5);

\node at (6.05,0.5) {{\tiny at $g=0$}};
\node at (14.05,0.5) {{\tiny at $g=0$}};
\node at (6.3,3) {{\tiny  at $g=\gamma \neq 0$}};
\node at (14.3,3) {{\tiny  at $g=\gamma \neq 0$}};

\node at (6.3,1) {{\footnotesize $\mathbb{A}^3_{(z_1,z_2,z_3)}$}};
\node at (14.3,1) {{\footnotesize $\mathbb{A}^3_{(z_1,z_2,z_3)}$}};
\node at (6.3,3.5) {{\footnotesize $\mathbb{A}^3_{(z_1,z_2,z_3)}$}};
\node at (14.3,3.5) {{\footnotesize $\mathbb{A}^3_{(z_1,z_2,z_3)}$}};

\node at (4.25,1.65) {{\tiny $Y_0$}};
\node at (12.25,1.65) {{\tiny $Y_0$}};
\node at (4.5,3) {{\tiny $Y_\gamma$}};
\node at (11.95,2.75) {{\tiny $Y_\gamma$}};

\node at (9.7,4.4) {{\tiny $C_1$}};
\node at (10.65,4.4) {{\tiny $C_2$}};
\node at (11.4,4.4) {{\tiny $C_3$}};

\draw[shift={(2.1,1.1)}, scale=0.006, draw opacity=.4, line width=\pgflinewidth+.9pt, smooth,variable=\t] plot ({15*\t*\t},{\t*\t*\t});

\draw[shift={(10.1,1.1)}, scale=0.006, draw opacity=.4, line width=\pgflinewidth+.9pt, smooth,variable=\t] plot ({15*\t*\t},{\t*\t*\t});

\draw[shift={(2.1,3.5)}, scale=0.1, draw opacity=.4, line width=\pgflinewidth+.9pt, rotate=270, smooth,variable=\t] plot ({\t},{\t*\t});

\draw[shift={(7.7,2.7)}, scale=1, draw opacity=.4, line width=\pgflinewidth+.9pt]
(4,0) 
to [out=90,in=190] (3,.5)
to [out=55,in=330](2.5,1)
to [out=0,in=300](3.5,1.5)
to [out=340,in=270](5,1.5);

\node (E) at (2.1, 1.1) {{\large $\bullet$}};
\node (F) at (10.1, 1.1) {{\large $\bullet$}};

\node (G) at (10.1, 3.73) {{\large $\bullet$}};
\node (H) at (11.2, 4.2) {{\large $\bullet$}};
\node (I) at (10.6, 3.12) {{\large $\bullet$}};

\path[draw,circle dotted,line width = .7mm] ( $  (9.9,4.45) $) to[out=-70,in=70] ++(-92:4.3);
\path[draw,circle dotted,line width = .7mm] ( $ (10.36,4.45) $) to[out=-70,in=60] ++(-101:4.3);
\path[draw,circle dotted,line width = .7mm] ( $ (11.15,4.44) $) to[out=280,in=95] ++(-105:4.4);
\end{tikzpicture}

 \vspace{.1in}

\begin{center}
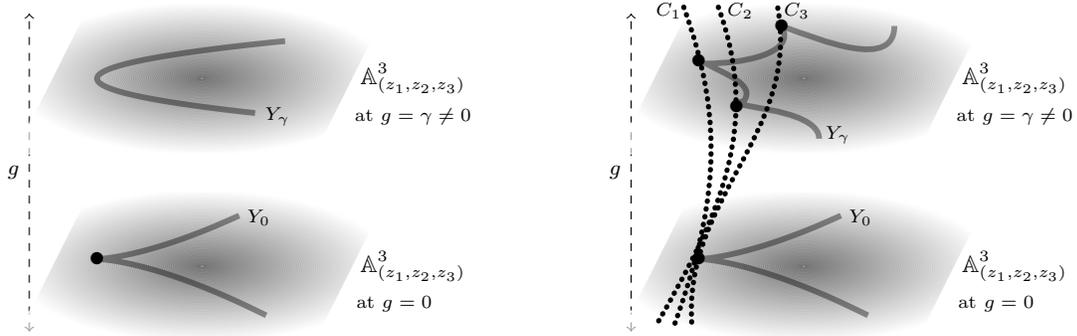
\textsc{Figure 1.} $Y$ for $3 \nmid n$ (left) and  $3|n$ (right). Bullets depict (curves  of) singularities.
\end{center}

\bigskip

\smallbreak Theorem~\ref{tintro2}(1) provides a strengthened version of \cite[Lemma~3.3]{BrownGoodearl}  without verifying the technical  {\it height 1 Azumaya hypothesis} of  \cite[Theorem~3.8]{BrownGoodearl}; it is also an extension of \cite[Theorem~1.3]{Walton} which was proved 
in the case when $(n,3)=1$.

\smallbreak Theorem~\ref{tintro2}(4) has the following immediate consequence:
\bco{remarkble}
Every 3-dimensional Sklyanin algebra that is module-finite over its center has 3 distinct isomorphism classes of central 
quotients if $(n,3) \neq 1$ and 2 of such isomorphism classes if~$(n,3)=1$. \qed
\eco

\green{The last result of the paper classifies the irreducible representations of the PI 3-dimensional Sklyanin algebras $S$ and their dimensions.
The irreducible representations of $S$ of intermediate dimension are those annihilated by $\mathfrak{m}_p \in \maxSpec Z$  corresponding to a point $p$ of $(C_1 \cup C_2 \cup C_3) \setminus \{\underline{0}\}$
in the case when the PI degree of $S$ is divisible by 3.}

\bth{interm} \green{If $(n,3)\neq 1$ and $p \in (C_1 \cup C_2 \cup C_3) \setminus \{\underline{0}\}$, then $S/(\mathfrak{m}_pS)$ has precisely 3 non-isomorphic  irreducible representations  and each of them has dimension $n/3$.}
\eth

\green{This fact was stated as a conjecture in the first version of the paper and a proof of it was given by Kevin DeLaet in \cite{DeLaet}. We provide an independent short proof in Section~\ref{proof-last-th}.}

%In the case when $(n,3) \neq 1$,  \cite[Theorem~1.3]{Walton} prompts the following conjecture.
%
%\bcj{interm} If $(n,3)\neq 1$ and $\mathfrak{m} \in (C_1 \cup C_2 \cup C_3) \setminus \{\underline{0}\}$, then $S/(\mathfrak{m}S)$ has precisely 3 non-isomorphic  irreducible representations  and each of them has dimension $n/3$.
%\ecj
%
%The resolution of the conjecture will fully describe the representation theory of the PI Sklyanin algebras. 
%We verify this in the appendix for an example of a 3-dimensional Sklyanin algebra $S$ of PI degree 6. We also refer the reader to recent work of Kevin DeLaet \cite{DeLaet} addressing this conjecture.

\bre{physics} 
The representation theoretic results of Theorem~\ref{tintro2} and Corollary~\ref{cremarkble}  contribute to String Theory, namely to the understanding of marginal supersymmetric deformations of the N = 4 super-Yang-Mills theory in four dimensions; see \cite{BJL}. The so-called {\it F-term constraints} on the moduli spaces of vacua of these deformations are given by representations of three-dimensional Sklyanin algebras $S$. In the case when $S$ is module-finite over its center, the irreducible representations of $S$ of maximum dimension correspond to {\it D-branes} in the {\it bulk} of the vacua, and  irreducible representations of smaller dimension correspond to {\it fractional D-branes} of the vacua.
\ere

\thref{intro2} also provides a key step towards the full description of the {\it discriminant ideals} of the PI Sklyanin algebras that are module-finite 
over their centers. This is described in Section~\ref{sec:future}.

\smallbreak The strategy of our  proof  of \thref{intro}(1) is to use {\it specialization of algebras} \cite[Section~2.2]{BrownGordon}
but with a new degree of flexibility: a simultaneous treatment of specializations of all possible {\it levels} $N$ (see Definition~\ref{dN-Pois-ord}). The previous work on the construction of Poisson orders by De Concini, Kac, Lyuba\-shenko, Procesi (for quantum groups at root of unity) \cite{DKP, DL, DP} and by Brown, Gordon (for symplectic reflection algebras) \cite{BrownGordon} always employed {\it first level} specialization and PBW bases. But our approach circumvents the problem that Sklyanin algebras are not easily handled with noncommutative Gr\"{o}bner/ PBW basis techniques. To proceed, we define a family of 
{\it formal Sklyanin algebras} $S_\hbar$ which specialize to the 3-dimensional Sklyanin algebras at $\hbar = 0$, and extend this to specializations 
for all components of the algebro-geometric description of $S$ via twisted homogeneous coordinate rings. Then we consider 
sections $\iota : Z(S) \to S_\hbar$ of the specialization $S_\hbar \twoheadrightarrow S$ such that $\iota(Z(S))$ is central in $S_\hbar$
modulo $\hbar^N S_\hbar$. This induces structures of Poisson $Z(S)$-orders on $S$ by `dividing by $\hbar^N$'. An analysis of the highest 
non-trivial level $N$ and a proof that such exists leads to the desired 
Poisson order structure in Theorem~\ref{tintro}(1).

\smallbreak The formula for the Poisson bracket on $Z(S)$ (\thref{intro}(2)) is then obtained by using the fact $S$ is a {\it maximal order} \cite{Stafford} and by showing that the singular locus of $Y$ has codimension $\geq 2$ in $Y$; thus we can `clear denominators' for computations in $Z(S)$ to obtain Theorem~\ref{tintro}(2). 

\smallbreak  We give a direct proof of the classification of symplectic cores in \thref{intro2}(2,3) for all fields $\kk$, based on Theorem~\ref{tintro}(2). 
The proof of \thref{intro2}(1,4) is carried out in two stages. In the first step, we prove these results in the case $\kk = \mathbb{C}$ 
by employing the aforementioned result of Brown and Gordon \cite{BrownGordon}; 
see Theorem~\ref{tBrGor}. The results in \cite{BrownGordon} rely on integration of Hamiltonian flows and 
essentially use the hypothesis that $\kk = \mathbb{C}$. We then establish \thref{intro2}(1,4) using base change arguments 
and general facts on the structure of  Azumaya loci. 

\smallbreak The paper is organized as follows. Background material and preliminary results on Poisson orders, symplectic cores, and 3-dimensional Sklyanin algebras (\green{namely, their} good basis and Heisenberg group symmetry) are provided in Section~\ref{sec:back}. We then establish the setting for the proof of Theorem~\ref{tintro}(1) in Section~\ref{sec:formal}. We prove Theorem~\ref{tintro}(1), Theorem~\ref{tintro}(2), and \green{Theorems~\ref{tintro2} and \ref{tinterm}} in Sections~\ref{sec:constr-Pord}, \ref{sec:bracket}, and~\ref{sec:repthy}, respectively.  Further directions and additional results are discussed in Section~\ref{sec:future}, including connections to noncommutative discriminants. Explicit examples illustrating the results above for cases  $n=2$ and $n=6$ are provided in an appendix (Section~\ref{sec:2and6}).
\medskip
\\
\noindent
{\bf Acknowledgements.} We are thankful to Kenny Brown and Paul Smith for helpful discussions.
\green{The authors would also like to thank the anonymous referee for their valuable and thorough feedback that greatly improved the exposition of this manuscript.}
%%%%%%%%%%%%%%%%%
%%%%%%%%%%%%%%%%%
%%%%%%%%%%%%%%%%%
%%%%%%%%%%%%%%%%%

\sectionnew{Background material and preliminary results} \label{sec:back}

We provide in this section background material and preliminary results on Poisson orders, including the process of specialization mentioned in the introduction, \green{as well as} material on symplectic cores and on (the noncommutative projective algebraic geometry of) 3-dimensional Sklyanin algebras.

%%%%%%%%%%%%%%%%%

\subsection{Poisson orders and specialization} \label{sec:P-order}

Here we collect some definitions and facts about Poisson orders and describe an extension of the specialization technique 
for obtaining such structures. 

\smallbreak Let $A$ be a $\kk$-algebra which is module-finite over a central subalgebra $C$. We will denote by $\Der(A/C)$ the algebra of $\kk$-derivations of $A$ that preserve $C$. The following definition is due to Brown and Gordon \cite{BrownGordon}:

\bde{Pois-order} 
\begin{enumerate}
\item The algebra $A$ is called a {\it Poisson $C$-order} if there exists a $\kk$-linear map
$\del : C \to \Der(A/C)$ such that the induced bracket $\{.,.\}$ on $C$ given~by 
\begin{equation}
\label{Poisson}
\{z, z' \}:= \partial_{z}(z'), \quad z, z' \in C
\end{equation}
makes $C$ a Poisson algebra. 
The triple $(A,C, \partial \colon C \to \Der(A/C))$ will be also called a {\it Poisson order} in places where the role of $\partial$ 
needs to be emphasized.
\smallskip

\item Assume further that  an algebraic group $G$ acts rationally by algebra automorphisms on $A$, so that $C$ is preserved under this action. We say that $A$ is a {\it $G$-equivariant Poisson $C$-order} if $\partial_{g \cdot z}(a) = g \cdot \partial_z(g^{-1}\cdot a)$ for all $g \in G$, $z \in C$, $a \in A$.
\end{enumerate}
\ede

As discussed in \cite[Section~2.2]{BrownGordon}, specializations of families of algebras give rise to Poisson orders. 
Below we generalize this construction to obtain Poisson orders from higher degree terms in the derivation $\partial$.

\smallbreak
Let $R$ be an algebra over $\kk$ and $\hbar$ be a central element of $R$ which is regular, i.e., not a zero-divisor of $R$. Let $[r, r']:= r r' - r' r$ for $r, r' \in R$.

\bde{spec'n}
We refer to the $\kk$-algebra $R_0 := R/\hbar R$ as the {\it specialization of $R$ at $\hbar \in Z(R)$}.
\ede

Let $\theta : R \twoheadrightarrow R_0$ be the canonical projection; so, $\ker \theta = \hbar R$.
Fix a linear map $\iota : Z(R_0) \hra R$ such that $\theta \circ \iota = \id_{Z(R_0)}$. Let $N \in \Zset_+$ be such that 
\begin{equation}
\label{N}
[\iota(z), r] \in \hbar^N R \quad \mbox{for all} \quad z \in Z(R_0), ~r \in R.
\end{equation}

Note that this condition holds for $N=1$: take $\wt{w} \in \theta^{-1}(w)$ for $w \in R_0$ and we get  $\theta([\iota(z), \wt{w}]) = [z, w] =0$; further,
$\ker \theta = \hbar R$.

\bde{spec-deriv} For $w \in R_0$ and $z \in Z(R_0)$, the {\it special derivation of level $N$} is defined as 
\begin{equation}
\label{def-special}
\partial_z (w) := \theta \left( \frac{ [\iota(z), \wt{w}]}{\hbar^N}\right), \quad \mbox{where} \; \; \wt{w} \in \theta^{-1}(w). 
\end{equation}
\ede

The above is well-defined: if $\wt{w}_1, \wt{w}_2 \in \theta^{-1}(w)$, then $\wt{w}_1 - \wt{w}_2 \in \ker \theta = \hbar R$; thus 
$w':=(\wt{w}_1 - \wt{w}_2)/\hbar \in R$ and the assumption \eqref{N} implies that  
\[
\theta \left( \frac{ [\iota(z), \wt{w}_1] }{\hbar^N}\right) - \theta \left( \frac{ [\iota(z), \wt{w}_2] }{\hbar^N}\right) = 
\theta \left( \frac{ [\iota(z), w'] }{\hbar^{N-1}}\right) = 0. 
\]
We will also show in the next result that $\partial_z$ is indeed a derivation.

\smallbreak The following result extends the specialization method of \cite{DKP,Hayashi} that produced Poisson orders with 1st level special derivations.  In applications, on the other hand,
it will be essential for us to consider more general $N$.

\bpr{special1} Let $R$ be a $\kk$-algebra and $\hbar \in Z(R)$ be a regular element. Assume that $\iota : R_0 := R/(\hbar R) \hra R$ 
is a linear section of the specialization map $\theta : R \twoheadrightarrow R_0$ such that \eqref{N} holds for some $N \in \Zset_+$. 
Assume that $R_0$ is module-finite over $Z(R_0)$. 
\begin{enumerate}
\item If, for all $z \in Z(R_0)$, $\partial_z$ is a special derivation of level $N$, then  $$(R_0, Z(R_0), \partial \colon Z(R_0) \to \Der(R_0/ Z(R_0)))$$ is a Poisson order.

\smallskip

\item This Poisson order has the property that $\partial : Z(R_0) \to \Der(R_0/ Z(R_0))$ is a homomorphism of Lie algebras, i.e., 
\[
\partial_{z} \partial_{z'} (w) - \partial_{z'} \partial_{z} (w) = \partial_{ \{ z, z' \} } (w)
\quad \mbox{for} \quad z, z' \in Z(R_0), ~w \in R_0.
\]
\end{enumerate}
\epr 

\begin{proof} (1) First, we verify that $\partial_z \in \Der(R_0)$ for $z \in Z(R_0)$. For $i=1,2$, choose $w_i \in R_0$ and $\wt{w}_i \in \theta^{-1}(w_i)$. 
Then $\wt{w}_1 \wt{w}_2 \in \theta^{-1}(w_1 w_2)$ and 
\begin{align*}
\partial_z(w_1 w_2) &= \theta \left( \frac{ [\iota(z), \wt{w}_1 \wt{w}_2] }{\hbar^N}\right)
~= \theta \left( \frac{ [\iota(z), \wt{w}_1]}{\hbar^N} \wt{w}_2 \right) + 
\theta \left( \wt{w}_1 \frac{ [\iota(z), \wt{w}_2] }{\hbar^N}\right)\\ &= \partial_z(w_1) w_2 + w_1 \partial_z(w_2).
\end{align*}
Next, we check that $\partial_{z} (z') \in Z(R_0)$ for all $z, z' \in Z(R_0)$. For $w \in R_0$ and $\wt{w} \in \theta^{-1}(w)$,  
\begin{align*}
[\partial_{z} (z'), w] &= \theta \left( \frac{[ [\iota(z), \iota(z')], \wt{w}] }{\hbar^N} \right)
~ = \theta \left( \frac{ [\iota(z), [\iota(z'), \wt{w}]] }{\hbar^N} \right) - \theta \left( \frac{[ \iota(z'),[\iota(z), \wt{w}]] }{\hbar^N} \right) 
\\&= \partial_{z} \left( \theta([\iota(z'), \wt{w}]) \right) - \partial_{z'} \left( \theta ([\iota(z), \wt{w}]) \right) = {\partial_{z}([z', w]) - \partial_{z'}([z, w])}=0.
\end{align*}
Finally, the fact that the bracket
\[
\{z, z' \} := \partial_{z} (z'), \quad z, z' \in Z(R_0)  
\]
satisfies the Jacobi identity follows from the verification of  part (2) below.

\smallskip

(2) The assumption \eqref{N} implies that
\[
[\iota(z), [\iota(z'), r]] \in \hbar^{2 N} R \quad \mbox{for all} \quad z, z' \in Z(R_0),~ r \in R
\]
since $[\iota(z'), r] \in \hbar^N R$ and $[\iota(z), [\iota(z'), r]/\hbar^N] \in \hbar^N R$. For $\wt{w} \in \theta^{-1}(w)$, we have
\begin{align*}
\partial_{z} \partial_{z'} (w) - \partial_{z'} \partial_{z} (w) &=
\theta \left(  \frac{[\iota(z), [\iota(z'), \wt{w}]]}{\hbar^{2N}} \right) - \theta \left(  \frac{[\iota(z'), [\iota(z), \wt{w}]]}{\hbar^{2N}} \right)
\\
&= \theta \left(  \frac{[[\iota(z), \iota(z')], \wt{w}]}{\hbar^{2N}} \right) = \partial_{ \{z, z'\}} (w).
\end{align*}

\vspace{-.25in}
\end{proof}

\bco{special2} 
\begin{enumerate}
\item If, in the setting of \prref{special1}, $C \subset Z(R_0)$ is a Poisson subalgebra of $Z(R_0)$ with respect to the Poisson structure \eqref{Poisson}
and $R_0$ is module-finite over $C$, then $R_0$ is a Poisson $C$-order via the restriction of $\partial$ to $C$. 
\smallskip

\item If, further, the restricted section $\iota : C \hra R$ is an algebra homomorphism, then
\[
\partial_{z z'}(w) = z \partial_{z'}(w) + z' \partial_{z}(w) \quad
\mbox{for} \quad z, z' \in C, w \in R_0.
\]
\end{enumerate}
\eco

\begin{proof} Part (1)  follows from \prref{special1}, and part (2) is straightforward to check.
\end{proof}

We end this part by assigning some terminology to the constructions above.

\bde{N-Pois-ord} The Poisson order produced in either Proposition~\ref{pspecial1} or Corollary~\ref{cspecial2} will be referred to as a {\it Poisson order of level $N$} when the level of the special derivation needs to be emphasized.
\ede

%%%%%%%%%

\subsection{Symplectic cores and the Brown-Gordon theorem}
\label{sec:BGor-thm}
Poisson orders can be used to establish isomorphisms for different central quotients 
of a PI algebra via a powerful theorem of Brown and Gordon \cite{BrownGordon}, which generalized previous results along these lines for quantum groups at roots of unity \cite{DKP,DP,DL}. The result relies on the notion of {\em{symplectic core}}, introduced in \cite{BrownGordon}. Let 
$(C, \{.,.\})$ be an affine Poisson algebra over $\kk$. For every ideal $I$ of $C$, there exists a unique maximal Poisson ideal 
$\PP(I)$ contained in $I$, and $\PP(I)$ is Poisson prime when $I$ is prime \cite[Lemma~6.2]{Goodearl}. Now we recall some terminology from \cite[Section 3.2]{BrownGordon}.

\bde{symplcore}[$\mathcal{P}(I)$, $\mathcal{C}(\mm)$] \quad (1) We refer to $\mathcal{P}(I)$ above as the {\it Poisson core} of $I$.
\begin{enumerate}
\item[(2)] We say that two maximal ideals $\mm, \nn \in \maxSpec C$ of an affine Poisson algebra $(C, \{.,.\})$ are {\it equivalent} if $\PP(\mm) = \PP(\nn)$. 
\smallskip
\item[(3)] The equivalence class of $\mm \in \maxSpec C$ is referred to as the  {\em{symplectic core}} of $\mm$, denoted by $\mathcal{C}(\mm)$. The corresponding partition of $\maxSpec C$ is called {\em{symplectic core partition}}. 
\end{enumerate}
\ede

Algebro-geometric properties of symplectic cores (that they are locally closed and smooth) are proved in \cite[Lemma 3.3]{BrownGordon}.
For instance, in the case when $\kk = \mathbb{C}$ and  $C$ is the coordinate ring on a smooth Poisson variety $V$,  each symplectic leaf of $V$ is contained in a single symplectic 
core of $V = \maxSpec C$, i.e., the symplectic core partition of $V$ is a coarsening of the symplectic foliation of $V$.
Furthermore, by \cite[Theorem 7.4]{Goodearl}, the symplectic core of a 
point $\mm \in V$ is the Zariski closure $\ol{L}$ of the symplectic leaf $L$ through $\mm$ minus the union of the Zariski closures of all symplectic leaves 
properly contained in~$\ol{L}$. 

\smallbreak One main benefit of using the symplectic core partition is the striking  result below.

\bth{BrGor} \cite[Theorem~4.2]{BrownGordon} Assume that $\kk = \mathbb{C}$ and that $A$ is a Poisson $C$-order which is an affine $\Cset$-algebra. If $\mm, \nn \in \maxSpec C$ are in the same symplectic core, 
then there is an isomorphism between the corresponding finite-dimensional $\mathbb{C}$-algebras
\[
A/(\mm A) \cong A /(\nn A).
\]

 \vspace{-.25in}
 
  \qed
\eth

%%%%%%%%%%%%%%%%%

\subsection{Three-dimensional Sklyanin algebras}
\label{3DSkly}

Here we recall the definition and properties of the 3-dimensional Sklyanin algebras and the corresponding twisted homogeneous coordinate rings.

\bde{Skly3}[$S$, $S(a,b,c)$] \cite{ATV1} Consider the subset of twelve points $$\mathfrak{D}:=\{\pcoor{1:0:0}, \pcoor{0:1:0}, \pcoor{0:0:1}\} \cup \green{\{\pcoor{a:b:c}~|~a^3= b^3= c^3\}}$$ in $\mathbb{P}^2$.
The {\it 3-dimensional Sklyanin algebras} $S:= S(a,b,c)$ are $\kk$-algebras generated by noncommuting variables $x, y, z$ of degree 1 subject to relations
\[
ayz+bzy+cx^2 ~=~ azx+bxz+cy^2 ~=~ axy+byx+cz^2 ~=~ 0,
\]
for $\pcoor{a:b:c} \in \mathbb{P}^2_\kk \setminus \mathfrak{D}$ where 
 $abc \neq 0$ and $(3abc)^3 \neq (a^3+b^3+c^3)^3$.
\ede

\smallbreak These algebras come equipped with geometric data that is used to establish many of their nice ring-theoretic, homological, and representation-theoretic properties. 
To start, recall that a {\it point module} for a connected graded algebra $A$ is a cyclic, graded left $A$-module with Hilbert series $(1-t)^{-1}$; these play the role of points in noncommutative projective algebraic geometry. The parameterization of $A$-point modules is referred to as the {\it point scheme} of $A$. 

\bdl{geomS}\textnormal{[$E$, $\phi$, $\sigma$]} \cite[Eq.~1.6, 1.7]{ATV1}  The point scheme of the 3-dimensional Sklyanin algebra $S=S(a,b,c)$ is given by the elliptic curve
\begin{equation} \label{E0}
E:= \mathbb{V}\left(\phi \right) \subseteq \mathbb{P}^2_{\pcoor{v_1:v_2:v_3}}, \text{ where } \phi = (abc) (v_1^3+v_2^3+v_3^3)-(a^3+b^3+c^3)v_1v_2v_3.
\end{equation}
If $\pcoor{1:-1:0}$ is the origin of $E$, then there is an automorphism $\sigma= \sigma_{abc}$ of $E$ given by translation by the point $\pcoor{a:b:c}$. Namely,
\begin{equation} \label{sigma0}
\sigma\pcoor{v_1:v_2:v_3} = [ac v_2^2 - b^2 v_1 v_3~:~bc v_1^2 - a^2v_2 v_3~:~ ab v_3^2  - c^2 v_1 v_2].
\end{equation}
The triple $\left(E,  ~\mathcal{O}_{\mathbb{P}^2}(1)|_E, ~\sigma\right)$ is referred to as  the geometric data of $S$.
\qed
\edl

Using this  data, we consider a noncommutative coordinate ring of $E$; its generators are sections of the invertible sheaf 
$\mathcal{O}_{\mathbb{P}^2}(1)|_E$ and its multiplication depends on the automorphism $\sigma$.

\bde{thcr} Given a projective scheme $X$, an invertible sheaf $\mathcal{L}$ on $X$, and an automorphism $\sigma$ of $X$,
the {\it twisted homogeneous coordinate ring} attached to this geometric data is a graded $\kk$-algebra 
$$B(X,  \mathcal{L}, \sigma) = \textstyle \bigoplus_{i \geq 0} B_i, \quad B_i := H^0(X, \mathcal{L}_i)$$
 with $\mathcal{L}_0 = \mathcal{O}_X$, $\mathcal{L}_1 = \mathcal{L}$, and $\mathcal{L}_i = \mathcal{L} \otimes \mathcal{L}^{\sigma} \otimes \cdots \otimes \mathcal{L}^{\sigma^{i-1}}$ for $i \geq 2$. The multiplication map $B_i \otimes B_j \to B_{i+j}$ is defined by $b_i \otimes b_j \mapsto b_i b_j^{\sigma^i}$ using $\mathcal{L}_i \otimes \mathcal{L}_j^{\sigma^i} = \mathcal{L}_{i+j}$.
 \ede
 
\bnota{BELsig} [$B$, $\mathcal{L}$] From now on, let $B$ be the twisted homogeneous coordinate ring attached to the geometric data $\left(E, \mathcal L, \sigma\right)$ from \dlref{geomS}, where $\mathcal L$ denotes $\mathcal{O}_{\mathbb{P}^2}(1)|_{E}$.
 \enota

It is often useful to employ the following embedding of $B$ in a skew-Laurent extension of the function field of $E$.

\ble{embedB} \cite{AST}
Given $\left(E, \mathcal L,\sigma\right)$   from \dlref{geomS}, extend $\sigma$ to an automorphism of the field  $\kk(E)$ of rational functions on $E$ by $\nu^{\sigma}(p) = \nu(\sigma^{-1} p)$ for $\nu \in \kk(E)$, $p \in E$. 
For any nonzero section $s$ of $\mathcal L$, that is, any degree 1 element of $B$, take $D$ to be the divisor of zeros of $s$, and let $V$ denote $H^0(E, \mathcal{O}_{E}(D)) \subset \kk(E)$. 

\smallbreak Then, the vector space isomorphism $\nu s \mapsto \nu t$ for $\nu \in V$ extends to an embedding of $B$ in $\kk(E)[t^{\pm 1}; \sigma]$. Here, $t\nu = \nu^{\sigma} t$ for $\nu \in \kk(E)$. \qed
\ele

Now the first step in obtaining nice properties of 3-dimensional Sklyanin algebras is to use the result below.
 
\ble{factor-g} \textnormal{[$g$]} \cite[before Theorem~6.6, Theorem~6.8]{ATV1} \cite[(10.17)]{AS} 
The degree 1 spaces of $S$ and of $B$ are equal. Moreover, there is a surjective map from $S$ to $B$, whose kernel is generated by the regular, central degree 3  element below:
\begin{equation}
\label{g-elem}
g:= c(c^3-b^3)y^3+b(c^3-a^3)yxz+a(b^3-c^3)xyz+c(a^3-c^3)x^3.
\end{equation}

\vspace{-.25in}

\qed
\ele

Many good ring-theoretic and homological properties of $S$ are obtained by lifting such properties from the factor $B$,  some of which are listed in the following result.

\bpr{hilb} \cite{AS, ATV1}  The 3-dimensional Sklyanin algebras $S$ are Noetherian domains of global dimension 3 that satisfy the Artin-Schelter Gorenstein condition. In particular, the algebras $S$ have  Hilbert series $(1-t)^{-3}$.
%\green{\textnormal{[[}(REMOVED)\textnormal{]]}}
\qed
\epr
%precisely when $[a:b:c] \in \mathbb{P}^2_{\kk} \setminus \mathfrak{D}$.

Moreover, the representation theory of both $S$ and $B$ depend on the geometric data $\left(E,  \mathcal L, \sigma \right)$; this will be discussed further in the next section. For now, we have:

\bpr{PI}  \cite[Lemma~8.5]{ATV1} \cite[Theorem~7.1]{ATV2} Both of the algebras  $S$ and $B=B(E, \mathcal L, \sigma) \cong S/gS$ are module-finite over their center if and only if the automorphism $\sigma$ has finite order.  \qed
\epr

Hence, one expects that both $S$ and $B$ have  a large center when $|\sigma|<\infty$. (In contrast, it is well-known that the center of $S$ is $\kk[g]$, when $|\sigma|=\infty$.) Indeed, we have the following results. Note that we follow closely the notation of \cite{SmithTate}.

\ble{centerB} \textnormal{[$E''$, $\Phi$]} \cite[Lemma~2.2]{AST} \cite[Corollary~2.8]{SmithTate}  Given the geometric data $\left(E, \mathcal L, \sigma \right)$ from 
\dlref{geomS}, suppose that \green{$|\sigma|=:n$ with $1<n<\infty$.} 
%\green{\textnormal{[[}(REMOVED)\textnormal{]]}}
%\begin{enumerate}
%\item If $n=1$, then $B$ is commutative (and equal to its own center).
%\smallskip
%\end{enumerate}
Now take $E'':=E/\langle \sigma \rangle$, with defining equation $\Phi$, so that $E\to E''$ is a cyclic \'{e}tale cover of degree~$n$. Recall \leref{embedB} and let $D''$ be the image of $D$ on $E''$ and let $V''$ denote $H^0(E'', \mathcal{O}_{E''}(D''))$. 
%  \green{\textnormal{[[}(REMOVED)\textnormal{]]}}

%\smallskip

%\begin{enumerate}
%\item[(2)] If $n\neq 1$, 

Then the center of $B$ is the intersection of $B$ with $\kk(E'')[t^{\pm n}]$, which is equal to $\kk[V'' t^n]$, and this is also a twisted homogeneous coordinate ring of $E''$ for the embedding of $E'' \subseteq \mathbb{P}^2$ for which $D''$ is the intersection divisor of $E''$ with a line. \qed
%\end{enumerate}
\ele

Central elements of $B$ lift to central elements of $S$ as described below. 
We will identify $S_1 \cong B_1$ via the canonical projection. 

\bde{good}[$s$] \cite{AST,SmithTate} Let $s$ be the value $n/(n,3)$. A section 
of $B_1 := H^0(E, \mathcal{L})$ is called {\em{good}} if its divisor of zeros 
is invariant under $\sigma^s$ and
consists of 3 distinct points whose orbits under the group $\lcor \sigma \rcor$ 
do not intersect.  A {\em{good basis}} of $B_1$ is a basis that consists of good elements so that the $s$-th powers of these elements \green{generate $B_n^{\langle \sigma \rangle}$} if $(n,3)=1$ or generate $(B^{\langle \sigma^3 \rangle})_{n/3}$ if $3|n$.
\ede

The order of $\sigma^{s}$ equals $(n,3)$. 
As mentioned at \cite[top of page~31]{SmithTate},  $\sigma^s$ fixes the class $[\mathcal L]$ in $\mathrm{Pic}\, E$. 

\bnota{autom-rho}[$\rho$] 
The automorphism $\sigma^s$ of $E$ induces an automorphism of $B_1$ via the identification $\mathcal L^{\sigma^s}\cong \mathcal L$. 
This automorphism of $B$ will be denoted by  $\rho$. 
\enota

By \cite[Lemma 3.4]{SmithTate}, there is a unique lifting of $\rho$ to an automorphism of $S$ 
and the central regular element $g$ is $\rho$-invariant.  

\smallbreak Now we turn our attention to the Heisenberg group symmetry of $S$ and of $B$. Recall that the Heisenberg group $H_3$ is the group of upper triangular $3 \times 3$-matrices with entries in $\Zset_3$ and 1's on the diagonal. It acts by graded automorphisms on $S$ in such a way that $S_1$ is the standard 3-dimensional representation of $H_3$. More concretely, the generators of $H_3$
act by 
\[
\rho_1: (x,y,z) \mt  (\zeta x,\zeta y,\zeta z), \quad \rho_2: (x,y,z)\mt (x,\zeta y, \zeta^2 z), \quad \rho_3: (x,y,z) \mt (y,z,x).
\]
where $\zeta$ is a primitive third root of unity.

\ble{actionrho}
We have that $\rho\in H_3$. Furthermore, in the case $(n,3)=3$, $\rho$ has three distinct eigenvalues and any good basis of 
$B_1$ consists of three eigenvectors (i.e., is unique up to rescaling of the basis elements).
\ele
\begin{proof}
The first statement is trivial in the case $(n,3)=1$ since $\rho= \id$. Now, suppose that $(n,3)=3$. Note that $\rho_2$ fixes $g$ and \green{induces a translation} $E$. 
More precisely, $\rho_2([v_1:v_2:v_3])=[v_1:\zeta v_2:\zeta^2 v_3]$ is the translation on $E$ by the point $[1:-\zeta:0]$ with respect to the origin $[1:-1:0]$. This implies that 
$\rho=\sigma^s$ commutes with $\rho_2$ on $E$ which are both translations with respect to the same origin. \red{Assuming that $v_1,v_2,v_3$ represent the three coordinate functions in $\kk(E)$, then $\rho\rho_2(v_i)=\rho_2\rho(v_i)$. Now consider the actions of $\rho$ and $\rho_2$ on $B_1$ and act on $x_1,x_2, x_3$. Hence we get $\rho\rho_2(x_i)=\lambda \rho_2\rho(x_i)$ for some scalar $\lambda$ because they induce the same action on $E$.} Taking into account that $\rho^3=\rho_2^3=\mathrm{id}$ we obtain that \red{$\lambda=\zeta^i$ for some $0\leq i\leq 2$ and $\rho\rho_2=\zeta^i\rho_2\rho$ on $B_1$}. Similarly, one shows that 
$\rho\rho_3=\zeta^j\rho_3\rho$ for some $0\leq j\leq 2$. A straightforward analysis using both the eigenspaces of $\rho_2$ and the explicit formula of $\rho_3$ yields that $\rho\in H_3$.

\smallbreak The fact that $\rho$ has three distinct eigenvalues is in the proof of \cite[Lemma 3.4(b)]{SmithTate}. The fact that the good bases are $\rho$-invariant is derived from \cite[paragraph after Definition on page 31]{SmithTate}.
\end{proof}

For the reader's convenience, we include the following computations; these \green{will be used} only in the Appendix.

\bre{uniquegoodbasis}
Assume that $(n,3)=3$. Since good bases are $\rho$-invariant and $\rho \in H_3$, their form in terms of the standard basis is as follows:
\[
\begin{array}{lllll}
 \{x,y,z\}, &&&\text{ if } \rho=\rho_1^{\pm 1} \rho_2, &\rho_1^{\pm 1} \rho_2^2; \\
\{x+y+z, &x+\zeta^2 y+\zeta z, &x+\zeta y+\zeta^2 z\}, &\text{ if } \rho=\rho_1^{\pm 1} \rho_3, &\rho_1^{\pm 1} \rho_3^2;\\
 \{x+y+\zeta^2z, &x+\zeta^2 y+z, &x+\zeta y+\zeta z\}, &\text{ if } \rho=\rho_1^{\pm 1} \rho_2^2\rho_3, &\rho_1^{\pm 1} \rho_2\rho_3^2;\\
  \{x+y+\zeta z, &x+\zeta y+ z, &x+\zeta^2 y+\zeta^2 z\}, &\text{ if } \rho=\rho_1^{\pm 1} \rho_2\rho_3, &\rho_1^{\pm 1} \rho_2^2\rho_3^2. 
\end{array}
\]
\ere

\ble{goodprelim}

\begin{itemize}
\item[(1)] If $(n,3)=1$, then for every element $\tau \in \lcor \rho_2 \rcor \times \lcor \rho_3 \rcor \subset H_3$ of order 3, 
there exists a good basis $\{x_1,x_2,x_3\}$ of $B_1$ which is cyclically permuted by $\tau$: $\tau(x_i) = x_{i+1}$, indices taken modulo 3.
\item[(2)] If $(n,3)=3$, then each good basis of $B_1$ can be rescaled so that the action of the Heisenberg group $H_3$ 
on $B_1$ takes on the standard form for the 3-dimensional irreducible representation of $H_3$.
\end{itemize}
\ele

\begin{proof}
(1) \red{Note that the $n$-th power map $f: B_1\to B_n^{\langle \sigma\rangle}$ given by $x\mapsto x^n$ is surjective by \cite[Lemma 5]{AST}. Moreover, by \cite[Proposition 2.6]{SmithTate}, $ B_n^{\langle \sigma\rangle}$ can be naturally identified with $B(E/\langle \sigma\rangle,\sigma^n,\mathcal L')_1$, where $\mathcal L'$ is the descent of $\mathcal L_n$ to $E/\langle \sigma\rangle$. One can easily check that $\tau$ induces a translation on $E$ by some $3$-torsion point $p$ (e.g., $p=[1:-\zeta:0]$ when $\tau=\rho_2$). Hence $\tau$ gives a translation on $E/\langle \sigma\rangle$ by the image of $p$ via the surjection $E\twoheadrightarrow E/\langle \sigma\rangle$. Since $3\nmid|\sigma|=n$, we get $\tau$ is nontrivial on $E/\langle \sigma\rangle$; and hence is nontrivial on $B_n^{\langle \sigma\rangle}$ satisfying $\tau^3=1$. Since $\tau$ has three distinct eigenvalues both on $B_1$ and $B_n^{\langle \sigma\rangle}$, the proper $\tau$-invariant subspaces of both spaces are  three hyperplanes and three lines.} 

\smallskip
\red{Let $U$ be the subset of $B_1$ consisting of nonzero elements whose divisors of zeros consist of $3$ distinct points with the property that their orbits under the group $\langle \sigma\rangle$ do not intersect. Clearly, $U$ is a Zariski open dense subset of $B_1$. So we can take some $w\in U$ avoiding those proper $\tau$-invariant subspaces in $B_1$ and the inverse images of those through $f: B_1\to B_n^{\langle \sigma\rangle}$. We claim that $\{w,w^\tau,w^{\tau^2}\}$ is a good basis. First of all, by definition, $w, w^\tau, w^{\tau^2}$ are all good elements. Secondly, they are linearly independent in $B_1$; otherwise, $\text{span}(w, w^\tau, w^{\tau^2})$ is a proper $\tau$-invariant subspace of $B_1$, which contradicts to our choice of $w$. Similarly, one can show that the $n$-th powers of $w, w^\tau, w^{\tau^2}$ are also linearly independent in $B_n^{\langle \sigma\rangle}$. This completes the proof of part (1).}

\smallbreak (2) This part follows from the facts that $\rho \in H_3$ and that a good basis consists of $\rho$-eigenvectors with distinct eigenvalues; see \cite[paragraph after Definition on page~31]{SmithTate}.
\end{proof}

\bnota{notation-x} [$x_1$, $x_2$, $x_3$, $\tau$] For the rest of this paper we fix a good basis $x_1, x_2, x_3$ of $B_1$ with the properties in \leref{goodprelim} and an element 
$\tau \in H_3$ such that
\begin{equation}
\label{tau}
\tau(x_i) = x_{i+1}, \quad \mbox{for $i=1,2,3$, indices taken modulo 3} 
\end{equation}
One can see that $g$ is fixed under $\tau$ by direct computation.
\enota

Next recording some results of Artin, Smith, and Tate, we have the following.

\bpr{centerS} \textnormal{[$Z$, $z_1$, $z_2$, $z_3$, $F$, $u_1$, $u_2$, $u_3$, $\ell$, $E'$, $f_3$]} \cite{AST} \cite[Theorems 3.7, 4.6,~4.8]{SmithTate}  The center $Z$ of $S$ is given as follows.

\begin{enumerate}
%\item We have that $n>1$ (since $abc \neq 0$).
%\green{[ REMOVED]]} 
%
%\smallskip

\item The center $Z$ is generated by three algebraically independent elements $z_1, z_2, z_3$ of degree $n$ along with~$g$ in \eqref{g-elem}, 
subject to a single relation $F$ of degree 3n. In fact, there is a choice of generators $z_i$ of the form
\[
z_i = x_i^n + \textstyle \sum_{1 \leq j < n/3} c_{ij} g^j x_i^{n-3j}
\]
where $\{x_1, x_2, x_3 \}$ is a good basis of $B_1$ and $c_{ij} \in \kk$. 

\smallskip

\item  If $n$ is divisible by 3 then there exist elements $u_1$, $u_2$, $u_3$ of degree $n/3$
that generate the center of the Veronese subalgebra $B^{(n/3)}$ of $B$,
so that $z_i = u_i^3$.

\smallskip

\item If $n$ is coprime to 3, then 
$$F = g^n + \Phi(z_1,z_2,z_3),$$ where $\Phi$ is the degree 3 homogeneous polynomial defining $E''\subset \mathbb P^2$  in Lemma~\ref{lcenterB}; here, $z_1, z_2, z_3$ 
are the $n$-th powers of the {\it good basis} $\{x_1, x_2, x_3\}$ of $B_1$. 
\smallskip

\item If $n$ is divisible by 3, then 
$$F = ~~g^n +  3\ell g^{2n/3} + 3\ell^2 g^{n/3} + \Phi(z_1,z_2,z_3).$$ 
 Here, $\Phi$ is as in Lemma~\ref{lcenterB},  $\ell$ is a linear form vanishing on the three images in $E''$ of the nine inflection points of $E':=E/\langle \sigma^3 \rangle$, and $f_3$ is the defining equation of $E'\subset \mathbb P^2$. Moreover, $f_3(u_1,u_2,u_3) + g^s = 0$ in $Z(S^{(3)})$.
\qed
\end{enumerate}
\epr

\ble{Z3good}
For a good basis $x_1, x_2, x_3$ of $B_1$ with the properties in \leref{goodprelim},
 the coefficients $c_{ij}$ in \green{Proposition~\ref{pcenterS}(1)} only depend on $j$, that is
\begin{equation}
\label{good-u}
z_i = x_i^n + \textstyle \sum_{1 \leq j < n/3} c_{j} g^j x_i^{n-3j}, \quad \text{ for some $c_i \in \Bbbk$.}
\end{equation}
Furthermore, if $(n,3)=3$, then the 9-element normal subgroup $N_3 \cong \Zset_3 \times \Zset_3$ of $H_3$, that rescales $x_1, x_2, x_3$, 
fixes the center $Z$ of $S$. The linear form $\ell$ in \green{Proposition \ref{pcenterS}(4)} is given by $\ell=\alpha(z_1+z_2+z_3)$ for some $\alpha\in \kk^\times$.
\ele
\begin{proof}
The first part follows at once from Lemma \ref{lgoodprelim} and \green{Proposition \ref{pcenterS}(1)}.

\smallbreak Let $(n,3) =3$ and $\green{N_3=\langle \rho_1,\rho_2\rangle}\subset H_3$. In this case the elements of $N_3$ rescale each of the good basis elements $x_1, x_2, x_3$ and fix $g$. 
This implies that each of these elements fixes $z_1, z_2, z_3$. The last statement follows from the fact that $\ell$  is fixed by $\tau$.
\end{proof}
%%%%%%%%%%%%%%%%%
\sectionnew{A specialization setting for Sklyanin algebras} 
\label{sec:formal}
%%%%%%%%%%%%%%%%%
The goal of this section is to construct a specialization setting for the Sklyanin algebras that is compatible with the 
geometric constructions reviewed in Section \ref{3DSkly}. The section sets up  
some of the notation that we will use throughout this work. Fix 
$\pcoor{a:b:c:\al:\be:\ga} \in \mathbb{P}^5_\kk$ such that
$\pcoor{a:b:c} \in \mathbb{P}^2_\kk \setminus \mathfrak{D}$ satisfies the conditions of Definition~\ref{dSkly3}.

\bhy{standhyp}
In the rest of this paper, $S:= S(a,b,c)$ will denote a 3-dimensional Sklyanin algebra that is module-finite over its center $Z:=Z(S)$, 
so that $|\sigma|=:n$ with $1 <n < \infty$. Moreover, $B ~(\cong S/gS)$ will be the corresponding twisted homogeneous coordinate ring.
\ehy

\noindent The reader may wish to view Figure~2 at this point for a preview of the setting.

\subsection{The first and second columns in Figure 2.}
\label{sec:1-2Col}
The goal here is to construct a degree 0 deformation $\ol{S_\hbar}$ of $S$ using a formal parameter $\hbar$. The specialization map for $S$ will 
be realized via a canonical projection $\theta_S : \ol{S_\hbar} \to S$ given by $\hbar \mapsto 0$. Here, $\ol{S_\hbar}$ will have the structure of a $\kk[\hbar]$-algebra; the beginning of this section is devoted to ensuring that the construction of $\ol{S_\hbar}$ is $\kk[\hbar]$-torsion free.

To begin, set
\begin{equation}
\label{wt-abc}
\wt{a}:= a + \al \hb, \quad
\wt{b} := b + \be \hb, \quad
\wt{c} := c + \ga \hb \quad
\in \kk[\hbar].
\end{equation}
It is easy to check that 
\[
\pcoor{\wt{a}:\wt{b}:\wt{c}} \notin \mathfrak{D} (\ol{\kk(\hb)}).
\]

\bde{ext-formalS}[$S_\hbar$, $\wh{S}_\hbar$, $\ol{S}_\hbar$] Consider the following formal versions of $S$.
\begin{enumerate}
\item Define the {\em{extended formal Sklyanin algebra}} $S_\hb$ to be $\kk[\hbar]$-algebra
\[
S_\hb = \frac{\kk[\hb] \lcor x, y, z \rcor}{
( \wt{a}yz+\wt{b}zy+\wt{c}x^2, ~~\wt{a}zx+\wt{b}xz+\wt{c}y^2, ~~\wt{a}xy+\wt{b}yx+\wt{c}z^2)} \cdot
\]

\smallskip

\item Denote by $\wh{S}_\hbar$ the Sklyanin algebra over $\ol{\kk(\hb)}$
with parameters $(\wt{a}, \wt{b}, \wt{c})$.

\smallskip

\item Define the {\em{formal Sklyanin algebra}} to be the $\kk[\hb]$-subalgebra $\ol{S}_\hb$ of $\wh{S}_\hb$ generated by $x,y,z$, that is,
\[
\ol{S}_\hb := \kk \lcor \hb, x, y, z \rcor \subset \wh{S}_\hb.
\] 
\end{enumerate}
\ede

We view $S_\hb$ as a graded $\kk[\hb]$-module by setting $\deg x = \deg y = \deg z =1$. Each \green{graded component} $S_\hb$ 
is a finitely generated $\kk[\hb]$-module. Since $\kk[\hb]$ is a principal ideal domain, we can decompose
\[
(S_\hb)_n := F_n \oplus T_n
\]
where $F_n$ is a free (finite rank) $\kk[\hb]$-\blue{sub}module \blue{(nonuniquely defined)} and $T_n$ is \blue{the} torsion $\kk[\hb]$-\blue{sub}module. For $n=0$,
$F_0 = \kk[\hb]$ and $T_0=0$.

\smallbreak The three algebras above are related as follows. First,
\begin{equation}
\label{S-Swh}
\wh{S}_\hbar \cong S_\hbar \otimes_{\kk[\hbar]} \ol{\kk(\hb)}.
\end{equation}
Denote by $\tau : S_\hb \to \wh{S}_\hb$ the corresponding homomorphism. It follows from \eqref{S-Swh} that
\[
\ker \tau =  \textstyle \bigoplus_{n \geq 0} T_n.
\]
Moreover, the algebra $\ol{S}_\hb$ is given by 
\begin{equation}
\label{free-Sbar}
\ol{S}_\hb = \Im \tau \blue{ ~\cong ~} \textstyle \bigoplus_{n\geq 0} \blue{S_n/T_n}.
\end{equation}
Thus, the formal Sklyanin algebra $\ol{S}_\hb$ is \blue{a factor} of the extended formal Sklyanin algebra $S_\hb$ by its $\kk[\hb]$-torsion part.
%\blue{\textnormal{[[}(REMOVED)\textnormal{]]}}

\smallbreak Now we show how the three-dimensional Sklyanin algebras are obtained via specialization. For each $d \in \kk$, we have the specialization map 
\[
\theta_d : S_\hbar \twoheadrightarrow S_{a + \al d, b + \be d, c + \ga d} \quad \mbox{given by} \quad \hb \mt d 
\]
whose kernel equals $\ker \theta_d = (\hb -d) S_\hbar$. Set $\theta:=\theta_0$. We have the following result.

\ble{T-n} 
\begin{enumerate} \item We get $\rank_{\kk[\hb]} F_n = \green{\dim S_n}$.

\smallskip

\item For all $d \in \kk$ such that
\begin{equation}
\label{d-assum}
\pcoor{(a + \al d):(b + \be d):(c + \ga d)} \notin \mathfrak{D},
\end{equation}
we have that $(\hb - d)$ is a regular element of $S_\hb$ and $T_n = (\hb -d) T_n$.  

\smallbreak 
\item The specialization map $\theta : S_\hb \twoheadrightarrow S$
factors through the map $\tau : S_\hb \twoheadrightarrow \ol{S}_\hb$. 

\end{enumerate}
\ele

\begin{proof} (1) The algebras $S$ and $\wh{S}_\hb$ have the same Hilbert series $(1-t)^{-3}$. Equation \eqref{S-Swh} implies 
that 
\[
(\wh{S}_\hb)_n = (\ol{S}_\hb)_n \otimes_{\kk[\hb]} \ol{\kk(\hb)} \cong F_n \otimes_{\kk[\hb]} \ol{\kk(\hb)}.
\]
So, $\rank_{\kk[\hb]} F_n = \dim (\wh{S}_\hb)_n = \dim S_n.$

\smallbreak (2) Set $\wt{a}_d:= a + \al d$, $\wt{b}_d:= b + \be d$ and $\wt{c}_d:= c + \ga d$. The assumption \eqref{d-assum} implies
that the algebras $S$ and $S(\wt{a}_d, \wt{b}_d, \wt{c}_d)$ have the same Hilbert series; see Proposition~\ref{philb}. The surjectivity of the specialization 
map $\theta_d$ gives that 
\[
S(\wt{a}_d,\wt{b}_d,\wt{c}_d)_n = (F_n/(\hb-d)F_n) \oplus (T_n/(\hb-d)T_n).
\]
By part (1), $\dim S(\wt{a}_d,\wt{b}_d,\wt{c}_d)_n  = \dim F_n/(\hb-d)F_n$, so, $(\hb-d) T_n = T_n$. 
Since $T_n$ is a finitely generated torsion $\kk[\hb]$-module, it is a finite-dimensional $\kk$-vector space and 
\[
\dim \ker (\hb -d)|_{T_n} = \dim T_n - \dim \Im (\hb -d)|_{T_n} = \dim T_n - \dim T_n =0.
\]
Hence, $\hb -d$ is a regular element of $S_\hb$. 

\smallbreak (3) This follows from part (2) and \eqref{free-Sbar}.
\end{proof}

\bnota{thetaS}[$\theta_S$] Denote by $\theta_S$ the corresponding {\it specialization map for the formal Sklyanin algebra $\ol{S}_\hbar$}, namely
\[
\theta_S : \ol{S}_\hb \to S \quad \mbox{given by} \quad \hb \mt 0.
\]
\enota

These maps form the two leftmost columns of the diagram in Figure~2 and the above results show the commutativity of the cells of the diagram 
between the first and the second column.

%%%%%%%%%%%%
\subsection{The third column in Figure 2.}
\label{sec:3Col}
Now we want to extend the results in the previous section to (the appropriate versions of) twisted homogeneous coordinate rings. 

\bnota{newg}[$\wt{g}$] Denote by $\wt{g}$ the elements of $\wh{S}_\hb$ given by \eqref{g-elem}
with $(a,b,c)$ replaced by $(\wt{a},\wt{b},\wt{c})$. 
%\blue{\textnormal{[[}(REMOVED, do not define an element of $\wt{g}$ of $S_\hb$)\textnormal{]]}}
\enota
%\blue{\textnormal{[[}(OLD LEMMA 3.10 REMOVED)\textnormal{]]}}
%\ble{centerSh} The element $\wt{g}$ is in the center of $S_\hb$. 
%\ele
%\begin{proof} Let $s \in S_\hb$. \leref{factor-g} implies that $\theta_d(\wt{g}s - s \wt{g}) =0$ for all $d \in \kk$ satisfying~\eqref{d-assum}. 
%Therefore, $\wt{g} s - s \wt{g} \in S_\hb$ is divisible by $\hb - d$ for all such $d \in \kk$. All but finitely many elements 
%$d \in \kk$ satisfy \eqref{d-assum} and, by \leref{T-n}(2), for all such $d$ we have that $\hb -d \in S_\hb$ are regular elements. 
%Since $S_\hb$ is  direct sum of finitely generated $\kk[\hb]$-modules, we get that $\wt{g} s-s \wt{g}=0$; that is $\wt{g} \in Z(S_\hb)$.
%\end{proof}

\blue{The central property of $g$ in Lemma \ref{lfactor-g} implies that}
\begin{equation}
\label{g-wt-Z}
\blue{\wt{g} \in Z(\wh{S}_\hb) \cap \ol{S}_\hb = Z(\ol{S}_\hb).}
\end{equation}
\bde{formalthcr}[$E_\hb$, $\mathcal{L}_\hb$, $\sigma_\hb$, $\wh{B}_\hb$, $\ol{B}_\hb$]
Denote by $E_\hb$ the elliptic curve over $\ol{\kk(\hb)}$, by $\mathcal{L}_\hb$ the invertible sheaf over $E_\hb$, and by $\sigma_\hb$ 
the automorphism of $E_\hb$
corresponding to $\wh{S}_\hb$ as in \dlref{geomS} with $(a,b,c)$ replaced by $(\wt{a}, \wt{b}, \wt{c}$) from \eqref{wt-abc}. Let
\[
\wh{B}_\hb:=B(E_\hb, \mathcal{L}_\hb, \sigma_\hb)
\]
be the corresponding twisted homogeneous coordinate ring. Its subalgebra 
\[
\ol{B}_\hb := \kk \lcor \hb,x,y,z \rcor \subset \wh{B}_\hb,
\]
generated by 
$\hb, x,y,z$, will be called {\em{formal twisted homogeneous coordinate ring}}.
\ede

By \leref{factor-g}, we have a surjective homomorphism 
\[
\Psi_\hb : \wh{S}_\hb \twoheadrightarrow \wh{B}_\hb \quad \mbox{with kernel} \quad
\ker \Psi_\hb = \wh{S}_\hb \wt{g}.
\]
Its restriction to $\ol{S}_\hb$ gives rise to the surjective homomorphism
\begin{equation}
\label{psi}
\psi_\hb : \ol{S}_\hb \twoheadrightarrow \ol{B}_\hb.
\end{equation}
We have 
\[
\wh{B}_\hb = \ol{B}_\hb \otimes_{\kk[\hb]} \ol{\kk(\hb)} \quad \mbox{and} \quad
\wh{S}_\hb = \ol{S}_\hb \otimes_{\kk[\hb]} \ol{\kk(\hb)}.
\]
The map $\Psi_\hb$ is the induced from $\psi_\hb$ map via tensoring $-\otimes_{\kk[\hb]} \ol{\kk(\hb)}$.

\ble{formalthcr} The kernel of the homomorphism $\psi_\hb : \ol{S}_\hb \twoheadrightarrow \ol{B}_\hb$ is given by
$\ker \psi_\hb = \ol{S}_\hb \wt{g}.$
\ele

\begin{proof} Clearly, $\ker \psi_\hb \supseteq \ol{S}_\hb \wt{g}$. Since $\ker \psi_\hb = \oplus_{n \geq 0} (\ker \psi_\hb)_n$ 
and each $(\ker \psi_\hb)_n$ is a finite rank torsion-free $\kk[\hb]$-module, $\ker \psi_\hb$ is a free $\kk[\hb]$-module.
At the same time 
\[
\ker \psi_\hb \otimes_{\kk[\hb]} \ol{\kk(\hb)} ~\cong ~(\ol{S}_\hb \wt{g}) \otimes_{\kk[\hb]} \ol{\kk(\hb)}~ \cong ~\wh{S}_\hb \wt{g}.
\]
which is only possible if 
$\ker \psi_\hb = \ol{S}_\hb \wt{g}.$

\vspace{-.2in}

\end{proof}

\ble{thetaB} The composition $\ol{S}_\hb \stackrel{\theta_S}{\twoheadrightarrow} S\twoheadrightarrow B$ 
factors through the map $\psi_\hb : \ol{S}_\hb \twoheadrightarrow \ol{B}_\hb$.
\ele
\begin{proof} This follows from the description of $\ker \psi_\hb$ in \leref{formalthcr}.
\end{proof}

\bde{thetaB}[$\theta_B$]
Let $\theta_B : \ol{B}_\hb \twoheadrightarrow B$ be the map induced by Lemma~\ref{lthetaB}, which we  call the {\em{specialization map for the formal twisted homogeneous coordinate ring}} $\ol{B}_\hb$.
\ede

This completes the construction of the maps in the third column of the diagram in Figure~2 and proves the commutativity of 
its cells between the second and third column.  

%%%%%%%%%%%%
\subsection{The \green{fourth column} in Figure 2.} 

Now we complete Figure~2 as follows.
\label{sec:4Col}
\bde{Ahbar} [$A_\hbar$, $e$, $f$]
Denote by $A_\hbar$ the subring of the function field $\ol{\kk(\hbar)}(E_\hbar)$ generated by 
$e:= v_2/v_1$, $f:=v_3 /v_1$, and $\hbar$. 
\ede
The generators satisfy the relation
$
\widetilde{a}\widetilde{b}\widetilde{c}(1+e^3+f^3)-(\widetilde{a}^3+\widetilde{b}^3+\widetilde{c}^3) ef =0.$ 
%\green{[[REMOVED]]}
%The quotient field of $A_\hbar$ is 
%$Q(A_\hbar) \cong \ol{\kk(\hbar)}(E_\hbar).$

\bde{Rhbar} [$R_\hbar$] Let $R_\hbar := (A_\hbar)_{(\hbar)}$  be the {\it integral form} of the field $\ol{\kk(\hbar)}(E_\hbar)$.
\ede

Using  \eqref{sigma0} with replacing $(a,b,c)$ by $(\wt{a}, \wt{b}, \wt{c})$, one sees that the automorphism $\sigma_\hbar \in \Aut \left( \ol{\kk(\hbar)}(E_\hbar) \right)$ restricts to an automorphism of $R_\hbar$, 
given by 
\[
\sigma_\hbar(e)=\frac{ \wt{b} \,  \wt{c}- \wt{a}^2\, ef}{ \wt{a} \, \wt{c} \, e^2- \wt{b}^2\, f} \quad \text{ and } \quad 
\sigma_\hbar (f)=\frac{ \wt{a} \, \wt{b} \, f^2- \wt{c}^2\, e}{ \wt{a} \, \wt{c} \, e^2- \wt{b}^2\, f}.
\]
Similar to Lemma~\ref{lembedB}, we have the canonical embeddings
\[
\ol{B}_\hbar \hra R_\hbar [t^{\pm 1}; \sigma_\hbar] \hra \ol{\kk(E_\hbar)}[ t^{\pm 1}; \sigma_\hbar].
\]
The ring $R_\hbar [t^{\pm 1}; \sigma_\hbar]$ is a graded localization of $\ol{B}_\hbar$ by an Ore set 
which does not intersect the kernel $\ker \theta_B = \hbar \ol{B}_\hbar$. Therefore the following map is well-defined.

\bde{thetaR}[$\theta_R$] Let $\theta_R : R_\hbar[t^{\pm 1}; \sigma_\hbar] \twoheadrightarrow \kk(E)[t^{\pm 1}, \sigma]$ be defined by $$\theta_R(e) =e, \quad \theta_R(f) =f,  \quad \theta_R(t)=t, \quad \theta_R(\hbar)=0,$$
which is the extension of $\theta_B$ via localization. We also denote by $\theta_R$ its restriction to the specialization map $R_\hbar \twoheadrightarrow \kk(E)$. These maps are referred to as the {\it specialization maps for the integral form of the formal twisted homogeneous coordinate $B$.}
\ede

The commutativity of the cells in Figure~2 between the third and forth column 
follows directly from the definitions of the maps in them.

\[
\xymatrix@C=2.1em@R=4em{
S_\hb\otimes_{\kk[\hb]}\overline{\kk(\hb)} \ar[r]^-{\cong} & \widehat{S}_\hb=\frac{\overline{\kk(\hb)}\left\langle x,y,z\right\rangle}{\left(\text{rel}(\tilde{a},\tilde{b},\tilde{c})\right)}\ar@{->>}[r]^-{\text{mod}(\wt{g})}& \widehat{B}_\hb=B(E_\hb,\mathcal L_\hb,\sigma_\hb) \ar@{^{(}->}[r]^-{\text{g.q.r.}} &  \overline{\kk(\hb)}(E_\hb)[\,t^{\pm 1};\sigma_\hb]\\
S_\hb=\frac{\kk[\hb]\left\langle x,y,z\right\rangle}{\left(\text{rel}(\tilde{a},\tilde{b},\tilde{c})\right)}\ar@{->>}[r]^-{\tau}\ar[u]\ar@{->>}[dr]^-{\theta} & \overline{S}_\hb=\kk\langle \hb,x,y,z\rangle  \ar@{->>}[r]^-{\text{mod}(\wt{g})}\ar@{^{(}->}[u]\ar@{->>}[d]^-{\theta_S}& \overline{B}_\hb=\kk\langle \hb,x,y,z\rangle \ar@{^{(}->}[r]\ar@{^{(}->}[u]\ar@{->>}[d]^-{\theta_B} & R_\hb[\,t^{\pm 1};\sigma_\hb]\ar@{^{(}->}[u]\ar@{->>}[d]_-{\theta_R}\\
& S = S(a,b,c) \ar@{->>}[r]^-{\text{mod}(g)} & B\ar@{^{(}->}[r]^-{\text{gr. quot. ring}}  & \kk(E)[\,t^{\pm 1};\sigma]   \\
  & Z\left(S\right) \ar@{->>}[r]^-{\text{mod}(g)}\ar@{^{(}->}[u]&Z(B) \ar[r]\ar@{^{(}->}[u]\ar@{^{(}->}[r]^-{\text{gr. quot. ring}} & \green{\kk(E)^{\sigma}[\,t^{\pm n}]}\ar@{^{(}->}[u] 
}
\]
\vspace{.1in}
\begin{center}
\textsc{Figure 2.} Specialization setting for Sklyanin algebras.\\
Integral forms, Poisson orders, and centers are respectively in the last three rows.
\end{center}

\vspace{.1in}

%%%%%%%%%%%%%%%%%
\sectionnew{Construction of non-trivial Poisson orders on Sklyanin algebras}
\label{sec:constr-Pord}

%%%%%%%%%%%%%%%%%
In this section we construct Poisson orders on all PI Sklyanin algebras $S$ for which the induced Poisson structures 
on $Z(S)$ are nontrivial. This gives a proof of Theorem~\ref{tintro}(1). We also construct nontrivial Poisson orders on 
the corresponding twisted homogeneous coordinate ring $B$ and the skew polynomial extension $\kk(E)[t^{\pm 1}; \sigma]$,
such that the three Poisson orders are compatible with each other.

\smallbreak We use the notation from the previous sections, especially standing Hypothesis~\ref{nstandhyp}.

%%%%%%%%%%%%%%
\subsection{Construction of orders with nontrivial Poisson brackets}
\label{Constr-Theorem}
Denote by
\[
X_n := \{ \pcoor{a:b:c} \in \mathbb{P}^2_\kk \setminus \mathfrak{D} \mid \sigma_{abc}^n=1 \}
\]
the parametrizing set for the Sklyanin algebras of PI degrees which divide $n$
(recall \deref{Skly3} for the notation $\mathfrak{D}$).
Throughout the section we will assume that, for the fixed $\pcoor{a:b:c} \in X_n$,
\begin{equation}
\label{albega-cond}
\pcoor{\al:\be:\ga} \in \mathbb{P}^2_\kk \quad \mbox{is such that} \quad \pcoor{a + d \al:b+ d \be: c+ d \ga} \notin 
X_n \sqcup \mathfrak{D} \; \; \mbox{for some} \; \; d \in \kk.
\end{equation}
This defines a Zariski open subset of $\mathbb{P}^2_\kk$ because $X_n \sqcup \mathfrak{D}$ is a closed proper subset 
of $\mathbb{P}^2_\kk$.

\bnota{xxu}[$x_i$, $\wt{x}_i$, $\wt{z}_i$]
Fix a good basis $x_1, x_2, x_3$ of $B_1$ as in \leref{goodprelim}. Throughout we will identify $B_1$ with $S_1$. 
Denote by $\wt{x}_i$ their preimages under the specialization map
$\theta_S : \ol{S}_\hbar \twoheadrightarrow S$ which are given by the same linear combinations of the generators $x,y,z$ of $\ol{S}_\hbar$ 
as are $x_i$ given in terms of the generators $x,y,z$ of $S$. Denote
\begin{equation}
\label{wt-u}
\wt{z}_i := \wt{x}_i^n + \textstyle \sum_{1 \leq j < n/3} c_{j} \wt{g}^j \wt{x}_i^{n-3j} \in \ol{S}_\hbar
\end{equation}
for the scalars $c_{j} \in \kk$ from \eqref{good-u}.
\enota

\bde{good-ect} A degree 0 section $\iota : Z \hra \ol{S}_\hbar$ of the specialization map $\theta_S : \ol{S}_\hbar \twoheadrightarrow S$ will be called {\it good} if 
\begin{enumerate}
\item $\iota(z_i) - \wt{z}_i \in \wt{g} \cdot \kk \lcor \wt{x}_i, \wt{g}, \hbar \rcor$ for the elements from \eqref{wt-u}, with the same noncommutative 
polynomials in three variables for $i =1,2,3$, and
\item $\iota(g) = \wt{g}$. 
\end{enumerate}
\ede

Now we define specialization in this context.

\bde{good-spec} We say that the specialization map $\theta_S:~\ol{S}_\hb~\to~S$ is a {\it good specialization of $S$ of level $N$} if there exists a good section 
$\iota : Z \hra \ol{S}_\hbar$  such that 
\begin{equation}
\label{NN}
[\iota(z), w] \in \hbar^N \ol{S}_\hbar \quad \mbox{for all} \quad z \in Z,~ w \in \ol{S}_\hbar. 
\end{equation}
\ede

Note that for every section $\iota : Z \hra \ol{S}_\hbar$ of $\theta_S$,
\[
[\iota(z), w] \in \hbar \ol{S}_\hbar \quad \mbox{for all} \quad z \in Z, ~w \in \ol{S}_\hbar. 
\]
Therefore, $N \geq 1$. Now we show that for a fixed value $n$, the levels $N$ of good specializations for $S$ of PI degree $n$ is bounded above.

\ble{exist} If $\pcoor{\al:\be:\ga} \in \mathbb{P}^2_\kk$ satisfies \eqref{albega-cond} and 
$N$ is a positive integer satisfying \eqref{NN} for a good section of $\theta_S$, then the levels $N$ of good specializations for $S$ have an upper bound. 
\ele

\begin{proof} 
First, denote by $A_\hbar$ and $R_\hbar :=(A_\hbar)_{(\hbar)}$ the rings defined analogously to the ones in 
Section \ref{sec:4Col} with dehomogenization performed with respect to 
$x_1 \in B_1$ not $v_1$. The condition \eqref{albega-cond} implies that 
the automorphism $\sigma_\hbar$ of $R_\hbar$ does not have order dividing~$n$.
So, $R_\hbar^{\sigma_\hbar^n} \subsetneq R_\hbar$. Moreover,  $\cap_{m \in \Zset_+} \hbar^m R_\hbar =0$. This implies that there exists a least positive integer $M$ such that 
\begin{equation}
\label{M}
\nu^{\sigma_\hbar^n} -\nu \notin \hbar^M R_\hbar \quad \mbox{for some} \quad \nu \in R_\hbar.
\end{equation}
We claim that $M \geq 2$ and $N<M$.

\smallbreak Since 
\[
\nu^{\sigma_\hbar^n} -\nu \in \hbar R_\hbar \quad \mbox{for all} \quad \nu \in R_\hbar
\]
one sees that  $M \geq 2$.

\smallbreak Towards the inequality $N<M$, assume that $\iota : Z \hra \ol{S}_\hbar$ is a good section of $\theta_S$, satisfying~\eqref{NN} for some positive integer $N$.
Recall~\eqref{psi}. We have 
\[
[\psi_\hb \iota(z), w] \in \hbar^N \ol{B}_\hbar \quad \mbox{for all} \quad z \in Z, ~w \in \ol{B}_\hbar. 
\]
Using that $R_\hbar[t^{\pm 1}; \sigma_\hbar]$ is a localization of $\ol{B}_\hbar$, where $\ol{B}_\hb \hookrightarrow R_\hb[t^{\pm 1}; \sigma_\hb]$ sending $\wt{x}_1$ to $t$, we obtain 
\begin{equation}
\label{zw-comm}
[\psi_\hb \iota(z), w] \in \hbar^N R_\hbar[t^{\pm 1}; \sigma_\hbar] \quad \mbox{for all} \quad z \in Z, ~w \in R_\hbar[t^{\pm 1}; \sigma_\hbar]. 
\end{equation}

Since $\ker \psi_\hb = \wt{g} \ol{S}_\hbar$ and $\iota$ is a good section, we get $\psi_\hb \iota(z_1) = t^n \in R_\hbar[t^{\pm 1}; \sigma_\hbar]$. Applying~\eqref{zw-comm} to $z = z_1$ and \green{$\nu \in R_\hbar$ gives}
\[
[t^n,\nu]=(\nu^{\sigma_\hbar^n}-\nu)t^n\in \hbar^N R_\hbar[t^{\pm 1}; \sigma_h] \quad \mbox{for all} \quad \nu \in R_\hbar,
\]
and thus,
\[
\nu^{\sigma_\hbar^n} -\nu \in \hbar^N R_\hbar \quad \mbox{for all} \quad \nu \in R_\hbar.
\]
Therefore $N < M$. 
\end{proof}

The following theorem provides a construction of a Poisson order with the non-vanishing property in \thref{intro}(1). Recall from the Introduction the definition and action of the group $\Sigma := \mathbb{Z}_3 \times \kk^\times$.

\bth{constr-Pord} Assume that $S$ is a Sklyanin algebra of PI degree $n$ so that
$\pcoor{\al:\be:\ga} \in \mathbb{P}^2_\kk$ satisfies \eqref{albega-cond}. Then the Poisson 
order $(S, Z, \partial: Z \to \Der(S/Z))$ of level~$N$,
constructed via good specialization of maximum level $N$, is $\Sigma$-equivariant and has the property that the induced Poisson structure on $Z$ is non-zero.
\eth

The theorem is proved in Section \ref{sec:proof-tPord}. 
%%%%%%
\subsection{Derivations of PI Sklyanin algebras} 
\label{sec:deri}
For an element $r$ of an algebra $R$, we will denote by $\ad_r$ the corresponding inner derivation of $R$; that is, 
\[
\text{ad}_r(s) = [r,s], \quad s \in R.
\]
We will need the following general fact on derivations of Sklyanin algebras which will be derived 
from the results of Artin-Schelter-Tate \cite{AST} for $(n,3)=1$ and of Smith-Tate \cite{SmithTate} in general.

\bpr{deriv} Let $S$ be a PI Sklyanin algebra and, by abusing notation, let $x \in S_1$ be a good element. If $\delta \in \Der S$ is such that
\begin{center}
\textnormal{(i)} $\delta|_Z =0$, \quad \quad \quad \quad 
\textnormal{(ii)}  $\delta(x) =0$, \quad \quad \quad \quad  and
\textnormal{(iii)} $\deg \delta =n$, 
\end{center}
then 
\[
\delta = c_1 g \ad_{x^{n-3}} + \cdots + c_m g^m \ad_{x^{n-3m}}
\]
for some non-negative integer $m<n/3$ and $c_1, \ldots, c_m \in \kk$. 
\epr

\begin{proof} Denote the canonical projection $\psi : S \twoheadrightarrow B$ with $\ker \psi = g S$. 
Since $\delta (g) =0$, $\psi \delta$ descends to a derivation of $B$ which, by abuse of notation, will be denoted by the same composition.
We extend $\psi \delta$ to a derivation of the graded quotient ring $\kk(E)[x^{\pm 1}; \sigma]$ of $B$.
It follows from (i) that $\psi \delta |_{\kk(E)^\sigma} =0$. Since $\kk(E)$ is a 
finite and separable extension of $\kk(E)^\sigma$, $\psi \delta|_{\kk(E)}=0$.
Indeed, if $a \in \kk(E)$ and $q(t)$ is its minimal polynomial over $\kk(E)^\sigma$, then 
\[
\big( \psi \delta(a) \big) q'(a) =0
\] 
because $\kk(E)^\sigma$ is in the center of $\kk(E)[x^{\pm 1}; \sigma]$. Since $q'(a) \neq 0$ and $\kk(E)[x^{\pm 1}; \sigma]$ is a domain, 
$\psi \delta(a) =0$. 

\smallbreak Finally, $\psi \delta(x)=0$ by (ii). Thus 
$\psi \delta =0$ as derivations on $\kk(E)[x^{\pm 1}; \sigma]$. So, $\delta(S)$ is \green{contained in} $\ker \psi = g S$. 

\smallbreak We define
\[
\delta_1:= g^{-1} \delta \in \Der(S). 
\]
Assumptions (i) and (ii) on $\delta$ imply that $\psi \delta_1$ descends to a derivation of $B$
(to be denoted in the same way) and that
$\deg (\psi \delta_1 ) = n-3$ and $\psi \delta_1 (x) =0.$
Applying \cite[Theorem~3.3, taking $r=1$]{SmithTate}, we obtain that there exists $c_1 \in \kk$ such that
$\psi \delta_1 = c_1 \ad_{x^{n-3}}.$
Therefore, 
\[
g^{-1} \delta - c_1 \ad_{x^{n-3}} = \delta_1 - c_1 \ad_{x^{n-3}} \in \Der S 
\quad \mbox{and} \quad
g^{-1} \Big( \delta - c_1 g \ad_{x^{n-3}} \Big) (S) \subseteq g S.
\]

Continuing this process, denote the derivation 
\[
\delta_2 := g^{-2} \big( \delta - c_1 g \ad_{x^{n-3}} \big) \in \Der(S).
\]
Similar to the composition $\psi \delta_1$, we obtain that $\psi \delta_2$ descends to a derivation of $B$
and that
$\deg (\psi \delta_2 ) = n-6$ and  $\psi \delta_2 (x) =0.$
By \cite[Theorem~3.3, taking $r=1$]{SmithTate}, there exists $c_2 \in \kk$ such that
$\psi \delta_2 = c_2 \ad_{x^{n-6}}$ and 
\[
g^{-2}\big( \delta - c_1 g \ad_{x^{n-3}} - c_2 g^2 \ad_{x^{n-6}} \big) \in \Der (S),
 \quad g^{-2}\big( \delta - c_1 g \ad_{x^{n-3}} - c_2 g^2 \ad_{x^{n-6}} \big)(S) \subseteq gS.
\] 

\smallbreak Let $m \in \Nset$ be the integer such that $m< n/3 \leq m+1$.
Repeating the above argument, produces
$c_1, \ldots, c_m \in \kk$
such that
\[
\delta_{m+1} := g^{-(m+1)} \big( \delta - c_1 g \ad_{x^{n-3}} - \cdots - c_m g^m \ad_{x^{n-3m}} \big) \in \Der (S).
\]
Since $\delta(g)=0$ by (i), $\psi \delta_{m+1}$ descends to a derivation of $B$ of degree $n- 3(m+1) \leq 0$. By \cite[Theorem~3.3, taking $r=1$]{SmithTate}, 
$\psi \delta_{m+1} =0$, so
$\delta_{m+1}(S) \subseteq g S.$

\smallbreak
Now using that $\deg \delta_{m+1} \leq 0$ gives
$\delta_{m+1}(S_1) \subseteq S_{\leq 1} \cap g S =0.$
Thus $\delta_{m+1}(S_1)=0$, which implies that $\delta_{m+1}=0$ since $S$ is generated in degree 1.
This completes the proof of the proposition.
\end{proof}

%%%%%%
\subsection{Proof of \thref{constr-Pord}}
\label{sec:proof-tPord}
It follows from \leref{Z3good} and the first condition in the definition of a good section, that for every good section 
$\iota : Z \hra \ol{S}_\hbar$, the Poisson order obtained by specialization
from $\iota$ is $\Sigma$-equivariant.

\smallbreak
Next we prove that the Poisson order obtained from a good specialization of maximum level $N$ has a nonzero Poisson structure on $Z$.
Assume the opposite; that is, the induced Poisson structure on $Z$ from the Poisson order $\partial$ vanishes.
Then we assert that 
\begin{center}
\textnormal{(i)} $\partial_{z_i} |_{Z} =0$, \quad \quad \quad
\textnormal{(ii)} $\partial_{z_i}(x_i) = 0$, \quad \quad  \quad and
\textnormal{(iii)} $\deg \partial_{z_i} =n$.
\end{center}
The first condition follows from the assumed vanishing of the Poisson structure. The third condition follows from~\eqref{def-special}. Now by Definition~\ref{dgood-ect}, 
$\iota(z_i) \in \kk \lcor \wt{x}_i, \wt{g}, \hbar \rcor$. It follows \blue{from \eqref{g-wt-Z}} that $[\iota(z_i), \wt{x}_i]=0$, and thus by~\eqref{def-special} we have
\begin{equation} \label{blah}
\partial_{z_i} (x_i) = \theta ([\iota(z_i), \wt{x}_i]/\hbar^N) =0.  
\end{equation}
This verifies the second condition above. 

\smallbreak Applying \prref{deriv} for \blue{$\delta= \partial_{z_1}$, gives that there exist 
$c_1, \ldots, c_m\in \kk$ and  $m < n/3$ such that
\[
\partial_{z_1} = c_1 g \ad_{x_1^{n-3}} + \cdots + c_m g^m \ad_{x_1^{n-3m}}.
\]
Therefore 
\begin{equation}
\label{first-bra}
[\iota(z_1) - \hbar^N( c_1 \wt{g} \, \wt{x}_1^{n-3} + \cdots + c_m \wt{g}^m \wt{x}_1^{n-3m}), w ] \in \hbar^{N+1} \ol{S}_\hbar \quad 
\end{equation}
for all $w \in \ol{S_\hbar}$. Recall that the Heisenberg group $H_3$ acts on $\ol{S}_\hbar$ by algebra automorphisms.
\green{The  
elements $\wt{x}_i \in \ol{S}_\hbar$ are defined  by the same linear combinations of the generators $x,y,z$ 
of $\ol{S}_\hbar$ as are those that define $x_i$ in terms of the generators $x,y,z$ of $S$.} This definition and 
Lemma \ref{lgoodprelim} imply that one of the elements of $H_3$ will act on 
$\ol{S}_\hbar$ as an automorphism that cyclically permutes $\wt{x}_1$, $\wt{x}_2$ and $\wt{x}_3$. This automorphism applied to 
\eqref{first-bra} gives
\[
[\iota(z_i) - \hbar^N(c_1  \wt{g} \, \wt{x}_i^{n-3} + \cdots + c_m \wt{g}^m \wt{x}_i^{n-3m}), w ] \in \hbar^{N+1} \ol{S}_\hbar \quad 
\]
\green{for $i=2,3$} and $w \in \ol{S}_\hbar$. Therefore,
\begin{equation}
\label{iota1}
\partial_{z_i} = c_1 g \ad_{x_i^{n-3}} + \cdots + c_m g^m \ad_{x_i^{n-3m}} 
\end{equation}
for all $i=1,2,3$.
}
Define a new section $\iota\spcheck : Z \hra \ol{S}_\hbar$ by first setting
\begin{equation}
\label{iota2}
\iota\spcheck(g) := \wt{g} \quad \text{ and } \quad 
\iota\spcheck(z_i) := \iota(z_i) - \hbar^N(  c_1 \wt{g} \, \wt{x}_i^{n-3} + \cdots + c_m \wt{g}^m \, \wt{x}_i^{n-3m}).
\end{equation}
Since $\iota\spcheck$ will not necessarily extend to an algebra homomorphism, choose a $\kk$-basis of~$Z$, namely $\mathcal{F}:=\{g^{l_0} z_1^{l_1} z_2^{l_2} z_3^{l_3} \mid (l_0, \ldots, l_3) \in L \}$ for some 
$L \subset \Nset^4$ with $g, z_1, z_2, z_3 \in \mathcal{F}$ and set
\begin{equation}
\label{iota3}
\iota\spcheck \left( g^{l_0} z_1^{l_1} z_2^{l_2} z_3^{l_3} \right) = 
\iota\spcheck(g)^{l_0} \iota\spcheck(z_1)^{l_1} \iota\spcheck(z_2)^{l_2} \iota\spcheck(z_3)^{l_3} \quad \mbox{for} \quad (l_0, \ldots, l_3) \in L.
\end{equation}
It is obvious that $\iota\spcheck : Z \hra \ol{S}_\hbar$ is a good section of $\theta_S$.

\smallbreak Combining \eqref{iota1} and \eqref{iota2} gives that 
\begin{equation}
\label{contra}
[\iota\spcheck(z), w]  \in \hbar^{N+1} \ol{S}_\hbar
\end{equation}
for $z = z_1,z_2, z_3, g$ and all $w \in \ol{S}_\hbar$. Indeed, $[\wt{g},w]=0$ and thus is divisible by $\hbar^{N+1}$. Further,
\[
\begin{array}{ll}
[\iota\spcheck(z_i), w] &= \hb^N\left(\frac{[\iota(z_i),w]}{\hbar^N} + c_1^i\wt{g}[w,\wt{x}_i^{n-3}] + \cdots + c_m^i\wt{g}^m[w,\wt{x}_i^{n-3m}]\right),
\end{array}
\]
and applying $\theta$ to the term in parentheses yields
$$\theta\left(\frac{[\iota(z_i),w]}{\hbar^N} + c_1^i\wt{g}[w,\wt{x}_i^{n-3}] + \cdots + c_m^i\wt{g}^m[w,\wt{x}_i^{n-3m}]\right) = \partial_{z_i}(\theta(w)) - \partial_{z_i}(\theta(w)) = 0.$$
Thus the term in the parenthesis is divisible by $\hbar$ and $[\iota\spcheck(z_i), w]$ is divisible $\hbar^{N+1}$. 

\smallbreak Now  \eqref{iota3} implies 
that \eqref{contra} is satisfied for all $z \in \mathcal{F}$ and $w \in \ol{S}_\hbar$, i.e., it holds for all $z \in Z$ and $w \in \ol{S}_\hbar$. 
This contradicts the hypothesis that $N$ is the maximum value satisfying~\eqref{NN}. Therefore, the induced Poisson structure on $Z$ from the Poisson order $\partial$ is nonvanishing.
\qed

%%%%%%%%%%%%%%%%%
\sectionnew{The structure of Poisson orders on Sklyanin algebras} \label{sec:bracket}
%%%%%%%%%%%%%%%%%
In this section, we establish \thref{intro}(2) on the induced Poisson bracket on the center $Z:=Z(S)$ of a 3-dimensional Sklyanin algebra $S$ via the specialization procedure described in Sections~\ref{sec:P-order} and~\ref{sec:1-2Col}. We make the following assumption for this section.

\bhy{goodhyp}
The generators $z_1, z_2, z_3$ of $Z$ are of the form \eqref{good-u}, i.e. given in terms of a good basis of $B_1 \cong S_1$.
\ehy

We begin with a straight-forward result.

\ble{g-Pcenter} For every good specialization of $S$ of level $N$, we have $\partial_g =0$ for the corresponding Poisson order 
on $S$. In particular, $g$ lies in the Poisson center of $Z$.
\ele

\begin{proof} 
A good section $\iota : Z \hra \ol{S}_\hbar$ has the property that $\iota(g) = \wt{g}$; see Definition~\ref{dgood-ect}. 
For $w \in S$ and $\wt{w} \in \theta_S^{-1} (w)$, \eqref{def-special} implies that
$\partial_g(w) = \theta_S\left([\wt{g},\wt{w}]/\hbar^N\right)$.
\green{By \eqref{g-wt-Z} we have that  $\wt{g} \in Z(\ol{S}_\hbar)$.} So, $\partial_g=0$. The last statement follows from \eqref{Poisson}.
\end{proof}

\subsection{The singular locus of $Y$}
\label{5.1}
Recall the notation from Section~\ref{3DSkly}. 
We  need the  preliminary result below.

\ble{pre-bracket} The following statements about $Y = \mathbb{V}(F)$ hold.
\begin{enumerate}
\item The partial derivatives $\partial_{z_i} F$ are nonzero for $i=1,2,3$.
\item The coordinate ring $\kk[Y]$ is equal to $\bigcap_{i=1}^3 \kk[Y][(\partial_{z_i} F)^{-1}]$.
\end{enumerate}
\ele

\begin{proof}
(1) This holds by a straight-forward computation using \green{Proposition~\ref{pcenterS}(3,4)}.

\smallbreak (2) Since $S$ is Noetherian, Auslander-regular, Cohen-Macaulay, and {\it stably free} by \cite{ATV1, Levasseur}, we can employ work of Stafford \cite{Stafford} to understand the structures of $Y$ and its coordinate ring $Z(S)$. Namely, $S$ is a maximal order, and thus $Z(S)$ is integrally closed and $Y$ is a normal affine variety, by \cite[Corollary on p.2]{Stafford}. 

\smallbreak Now by \cite[Discussion after proof of Corollary~11.4]{Eisenbud}, it suffices to verify that 
$$Y' := \mathbb{V}(\partial_{z_1}F, ~\partial_{z_2}F, ~\partial_{z_3}F) \cap Y$$
has codimension $\geq 2$ in $Y$. The definition of $Y'$ implies that 
it is the union of the singular loci of the slices $Y_\gamma:= Y \cap \mathbb{V}(g - \gamma)$,
\[
Y' = \textstyle \coprod_{\gamma \in \kk} (Y_\gamma)^{\textnormal{sing}}.
\]
In the following lemma we explicitly describe $Y'$ and prove that it coincides with the 
singular locus $Y^{\textnormal{sing}}$ of $Y$. In fact, the description of $Y^{\textnormal{sing}} = Y'$ 
below implies that the dimension of $Y^{\textnormal{sing}}$ is $\leq 1$, 
and hence has codimension $\geq$ dim($Y$)$ - 1$. Thus, the codimension of $Y^{\textnormal{sing}}$ in $Y$ is $\geq 2$, as needed. 
\end{proof}

\ble{Ysing}
The singular locus $Y^{\textnormal{sing}}$ of $Y$ is the origin if $(n,3) = 1$, and is the union of three dilation-invariant 
curves intersecting the coordinate hyperplanes at the origin if $(n,3) \neq 1$; the dilation is given by \eqref{dilation}.
Furthermore, $Y^{\textnormal{sing}} = Y'$, that is
\[
(Y_\gamma)^{\textnormal{sing}} = Y_\gamma \cap Y^{\textnormal{sing}} \quad \mbox{for} \quad \gamma \in \kk.
\]
\ele

\begin{proof} 
\green{To verify that $Y^{\textnormal{sing}} = Y'$, we need to show that $F=0$ and $\{\partial_{z_i} F =0\}_{i=1,2,3}$  imply that $\partial_g F =0$. Indices are taken modulo~3 below.}

\smallbreak \green{Say $(n,3) =1$. Recall from Proposition~\ref{pcenterS}(3) that  $F=g^n + \Phi(z_1, z_2, z_3)$. So, $\partial_g F =0$ if and only if $g=0$. Since $\Phi$ is a homogeneous equation in the $z_i$ of degree~3, we have  $\sum_{i=1}^3 z_i \partial_{z_i}F = 3 \Phi$. Assuming that  $\{\partial_{z_i} F =0\}_{i=1,2,3}$, we get $\Phi = 0$. Further, assuming that $F=0$ yields $g=0$, as desired. 
}

\smallbreak \green{Take $(n,3) \neq 1$. Recall  that  $F=g^n +3\ell g^{2n/3} + 3\ell^2g^{n/3}+ \Phi(z_1, z_2, z_3)$ [Proposition~\ref{pcenterS}(4)]. According to the proof of \cite[Theorem~4.8]{SmithTate}, the element $\Phi$ of Proposition~\ref{pcenterS} is of the form $\ell^3 - \mu z_1 z_2 z_3$, for $\mu \in \kk^{\times}$ and $\ell =\alpha(z_1+z_2+z_3)$ a (nonzero) degree one homogeneous polynomial by Lemma \ref{lZ3good}.}
 
 \smallbreak \green{Note that  $\partial_g F =0$ is equivalent to $\ell = -g^{n/3}$, when $n \neq 3$. (The fact that $S$ is a domain along with Proposition~\ref{pcenterS}(4) is used for the latter equivalence.) Now $F=0$ implies that $(\ell+g^{n/3})^3 = \mu z_1z_2z_3$. So if $F=0$ and $\partial_{z_i} F =0$ for $i=1,2,3$, then 
$$\textstyle \sum_{i=1}^3 z_i \partial_{z_i}F = 3\ell^3 - 3(\ell+g^{n/3})^3 +6\ell^2g^{n/3} + 3\ell g^{2n/3} = -3g^{n/3}(\ell+g^{n/3})^2 = 0.$$
Since $S$ is a domain, this yields $\ell = -g^{n/3}$ in the case $n \neq 3$, as desired. }

 \smallbreak \green{Consider the case $n=3$. If $F = 0$, then $\mu z_1 z_2 z_3 = (g+\ell)^3$. Moreover, if $\partial_{z_i}F = 0$, then $\mu z_{i+1} z_{i+2} = 3 \alpha (g + \ell)^2$ (indices taken modulo 3). So if $F = 0$ and $\partial_{z_i}F  = 0$ for all~$i$, then we obtain that $g + \ell = 0$, and this implies $\partial_{g}F  =0$.}

\smallbreak \green{Thus, $Y^{\textnormal{sing}} = Y'$. When $(n,3) =1$ we get $Y^{\text{sing}}  = Y^{\text{sing}}_0 = \mathbb{V}(\Phi) \cap \mathbb{V}(g)$, which is the origin since $E''$ is smooth; see Proposition~\ref{pcenterS}(3). In general, we make use of the fact that $E'' = \mathbb{V}(\Phi)$ is a smooth projective variety (thus, excluding a finite set of pairs $(\alpha, \mu)$ from calculations) to yield:
$$
Y^{\textnormal{sing}} = 
\begin{cases}
\{(0,0,0,0)\}, &\quad \text{ if } (n,3)=1\\
\bigcup_{i=1}^3 \{g^{n/3} + \alpha z_i = 0, ~~z_{i+1} = z_{i+2} =0\}, & \quad \text{ if } (n,3)\neq1.\\
\end{cases}
$$}

\vspace{-.25in}
\end{proof}

\subsection{Poisson structures on $Z(S)$}
\label{5.2}
Now we establish the main result of this section.

\bpr{bracket}
\begin{enumerate} 
\item Each homogeneous Poisson bracket on $Z(S)$, such that  $g$ is
in the Poisson center of $Z(S)$, is given by
\begin{equation}
\label{br}
\{z_1,z_2\} = \eta \partial_{z_3} F, \quad \{z_2,z_3\} = \eta \partial_{z_1} F, \quad \{z_3,z_1\} = \eta \partial_{z_2} F, \quad \text{ for some } \eta \in \kk.
\end{equation}
\item In the case when $S$ arises as a Poisson order of level $N$ (via Proposition~\ref{pspecial1} and Section~\ref{sec:1-2Col}), we get the formula above with $\eta \neq 0$, and obtain Theorem~\ref{tintro}(2) by rescaling.
\end{enumerate}
\epr

\begin{proof} 
(1) By \prref{centerS}(4,5) and the Leibniz rule, 
any Poisson structure on $Z(S)$ with $g$ in the Poisson center satisfies $\{F, z_i\} = \sum_{j=1}^3 (\partial_{z_j} F)\{z_j, z_i\}$; namely, $\{-, z_i\}$ is a derivation. So the following equations hold because $F = 0$ in $Z(S)$:
\[
\begin{array}{lll}
0 &= \{F,z_1\} &= (\partial_{z_2} F) \{z_2,z_1\} + (\partial_{z_3} F) \{z_3, z_1\},\\
0 &= \{F,z_2\} &= (\partial_{z_1} F) \{z_1,z_2\} + (\partial_{z_3} F) \{z_3, z_2\},\\
0 &= \{F,z_3\} &= (\partial_{z_1} F) \{z_1,z_3\} + (\partial_{z_2} F) \{z_2, z_3\}.
\end{array}
\]
The second equation with Lemma~\ref{lpre-bracket}(1) implies that $\{z_1,z_2\} = \eta (\partial_{z_3} F)$ for $\eta \in \kk(Y)_0$, since the bracket is homogeneous, $\deg(z_i)=n$, the left-hand side has degree $2n$,  and $\deg(F)=3n$. Now the first equation and Lemma~\ref{lpre-bracket}(1) imply that 
$$\{z_3,z_1\} = \frac{- (\partial_{z_2} F)\{z_2,z_1\}}{\partial_{z_3}F} = \eta (\partial_{z_2} F).$$
Likewise, $\{z_2,z_3\} = \eta (\partial_{z_1} F)$. Lemma~\ref{lpre-bracket}(2) allows us to clear denominators to conclude that $\eta \in \kk[Y]$. Therefore, $\eta \in  \kk(Y)_0 \cap \kk[Y] = \kk[Y]_0$, so $\eta \in \kk$. 

\smallbreak (2) Such an induced Poisson structure on $Z(S)$ is homogeneous since the formal Sklyanin algebra $\ol{S}_\hb$ is graded and since the bracket is given by \eqref{def-special}. Moreover, $g$ is in the Poisson center in this case by Lemma~\ref{lg-Pcenter}. So, this part follows from part~(1) and Theorem~\ref{tconstr-Pord}.
\end{proof}

\bre{otherPoisson} One can adjust the specialization method from Section~\ref{sec:formal} to involve different types of deformations of 
3-dimensional Sklyanin algebras $S$, such as the PBW deformations that appear in work of Cassidy-Shelton \cite{CassidyShelton}; 
see also work of Le Bruyn-Smith-van den Bergh \cite{LSV}. Unlike our setting above where the deformation parameter $\hbar$ has degree 0, 
the deformation parameter in the aforementioned works have either degree 1 or 2, which yields a Poisson bracket on $Z(S)$ of degree $-1$ or $-2$, 
respectively. This is worth further investigation.
\ere
%%%%%%%%%%%%%%%%%
%%%%%%%%%%%%%%%%%
%%%%%%%%%%%%%%%%%
%%%%%%%%%%%%%%%%%

\sectionnew{On the representation theory of $S$} \label{sec:repthy}
In this section we prove \thref{intro2}. Denote by $Y^1,Y^2, Y^3$, and $Y^4$ 
the strata of the partition of $Y = \maxSpec(Z(S))$ in \thref{intro2}(3). The stratum $Y^2$ is nonempty only in the case \blue{$(n,3) \neq 1$.}
Recall from the introduction that $\AA \subseteq Y$ denotes the Azumaya locus of $S$, and recall the group $\Sigma := \mathbb{Z}_3 \times \kk^{\times}$ acts on both $Z(S)$ and $Y$. 
\green{In the last part of the section we prove Theorem~\ref{tinterm}.}
%We end the section by providing a quick reduction of Conjecture~\ref{jinterm}.

\subsection{Symplectic cores and $\Sigma$-orbits}
\label{6.1} 
For $y \in Y$ denote by $\mm_y$ the corresponding maximal ideal of $Z = \kk[Y]$. 
Denote by $\CC(y)$ the symplectic core containing $y$.

\bpr{scores} Consider the Poisson structure on $Z$ given by 
\begin{equation}
\label{brackets}
\{z_1,z_2\} = \partial_{z_3} F, \quad \{z_2,z_3\} = \partial_{z_1} F, \quad \{z_3,z_1\} = \partial_{z_2} F,
\end{equation}
and $g$ is in the Poisson center of $Z$. Then, $Y_\gamma$ is a Poisson subvariety of $Y$ for each $\gamma \in \kk$. The symplectic cores of $Y$
are the sets $\{Y_\gamma \backslash Y^{\textnormal{sing}}_\gamma\}_{\gamma \in \kk}$, along with
the points in the union of curves $C_1 \cup C_2 \cup C_3$ if $(n,3) \neq 1$ or the point $\{\underline{0}\}$ if $(n,3)=1$.
\epr
\begin{proof} 
Since $g$ is in the Poisson center of $Z$, $(g - \gamma)$ is a Poisson ideal of $Z$, and thus $Y_\gamma$ is a Poisson subvariety of $Y$. 

\smallskip
The Poisson core $\PP(\mm_y)$ is a Poisson prime ideal by \cite[Lemma 6.2]{Goodearl}.
One easily sees that $y \in Y$ forms a singleton \blue{symplectic} core if and only if $\mm_y$ is a Poisson ideal, i.e., 
if all right hand sides of Poisson brackets in \eqref{brackets} belong to $\mm_y$. Lemma~\ref{lYsing} then implies that the singleton \blue{symplectic} cores are 
the points in the union of curves $C_1 \cup C_2 \cup C_3$ if $(n,3) \neq 1$ and the point $\{\underline{0}\}$ if $(n,3)=1$.

\smallskip
Let $y \in Y_\gamma \backslash Y^{\textnormal{sing}}_\gamma$ for some $\gamma \in \kk$. We have $\{y\} \subsetneq \mathbb{V}(\PP(\mm_y)) \subseteq Y_\gamma = \mathbb{V}(g - \gamma)$. Because $\PP(\mm_y)$ is a Poisson prime ideal, 
$\mathbb{V}(\PP(\mm_y)) \cap (Y_\gamma \backslash Y^{\textnormal{sing}}_\gamma)$ 
is a Poisson subvariety of the 2-dimensional smooth irreducible Poisson variety $Y_\gamma \backslash Y^{\textnormal{sing}}_\gamma$
whose Poisson structure is nowhere vanishing (i.e. it is a symplectic one). This is only possible if 
$\mathbb{V}(\PP(\mm_y))$ contains $Y_\gamma \backslash Y^{\textnormal{sing}}_\gamma$.
Therefore $\mathbb{V}(\PP(\mm_y))= Y_\gamma$, and thus $\PP(\mm_y) = (g-\gamma)$. This implies that 
$\CC(y) = Y_\gamma \backslash Y^{\textnormal{sing}}_\gamma$ which completes the proof of the proposition.
\end{proof} 

%%%%%%%%%%%%%%%%%%
\subsection{Proof of \thref{intro2}(2,3).} Part (2) follows from Lemma~\ref{lYsing} and Proposition~\ref{pscores}. Part (3) follows from \prref{scores}
and from the fact that with respect to the dilation action~\eqref{dilation}, we get
\[
\kk^\times \cdot (Y_\gamma \backslash Y^{\textnormal{sing}}_\gamma) = 
\begin{cases}
Y \backslash Y_0, &\mbox{if} \quad (n,3)=1
\\
Y \backslash (Y_0 \cup C_1 \cup C_2 \cup C_3), & \mbox{if} \quad (n,3) \neq 1
\end{cases}
\]
for all $\gamma \in \kk^\times$. Since the $\mathbb{Z}_3$-action cyclically permutes $C_1, C_2$ and $C_3$, 
$(C_1 \cup C_2 \cup C_3) \backslash \{ \underline{0} \}$ is a single $\Sigma$-orbit of symplectic 
cores.
\qed
%%%%%%%%%%%%%%%%%%
\subsection{Proof of \thref{intro2}(1,4) for $\kk = \Cset$}
\label{6.2} 
Using \thref{constr-Pord} and Proposition~\ref{pbracket}, we construct a Poisson order on $S$ for which the induced 
Poisson bracket on $Z$ is given by~\eqref{br} with $\eta \neq 0$. \thref{BrGor} and the fact that $\Sigma$ 
acts on $S$ by algebra automorphisms imply that 
\begin{equation}
\label{isomm}
S/(\mm_y S) \cong S/(\mm_{y'} S) \quad \mbox{if} \quad y,y' \in Y^j \text{ for some } j =1,2,3,4.
\end{equation}
This proves \thref{intro2}(4).

\smallbreak The stratum $Y^1$ and the Azumaya locus $\AA$ of $S$ are dense subsets of $Y$. Hence, $Y^1 \cap \AA \neq \varnothing$, 
and the isomorphisms
\eqref{isomm} imply that $Y^1 \subseteq \AA$. Thus, 
\[
S/(\mm_y S)  \cong M_n (\Cset) \quad \mbox{for} \quad y \in Y^1.
\]

The stratum $Y^3$ is a dense subset of \green{$Y_0 = \maxSpec (Z(B))$}. Since the Azumaya locus $\AA(B)$ is also 
dense in $Y_0$, the isomorphisms \eqref{isomm} imply that $Y^3 \subseteq \AA(B)$. 
The PI degree of $B = S/(g S)$ equals $n$ because the graded quotient ring of $B$ 
is isomorphic to $\Cset(E)[t^{\pm 1}; \sigma]$ and $\sigma$ has order $n$. 
Therefore $Y^3 \subseteq \AA(B) \subseteq \AA$ and 
 \[
S/(\mm_y S) \cong M_n (\Cset) \quad \mbox{for} \quad y \in Y^3.
\]

Finally, $Y^{\textnormal{sing}} \cap \AA = \varnothing$ by \cite[Lemma 3.3]{BrownGoodearl} and $Y^{\textnormal{sing}} = Y^2 \sqcup Y^4$ by Lemma~\ref{lYsing}.
Hence, 
\[
\AA = Y^{\textnormal{smooth}} = Y^1 \sqcup Y^3
\]
which proves \thref{intro2}(1).  \qed

%%%%%%%%%%%%%%%%%%

\subsection{Proof of \thref{intro2}(1,4)  for  $\kk = \ol{\kk}$ of characteristic 0}
\label{6.3} 
The set $Y^2$ is a single $\Sigma$-orbit, so \eqref{isomm} holds for $j=2$. The set $Y^4$ is a singleton.
The corresponding factor $S/(\mm_{\underline{0}} S)$ of $S$ is a finite dimensional algebra which is connected 
graded, and thus local. We also have $Y^{\textnormal{sing}} \cap \AA = \varnothing$ by \cite[Lemma 3.3]{BrownGoodearl}.

\smallbreak It suffices to show that for $y:=(y_1,y_2,y_3,y_4) \in Y \subset \mathbb{A}^4_{(z_1, z_2, z_3, g)}$ we have
\begin{equation} 
\label{factors}
S/(\mm_y S) \cong M_n(\kk) \quad \mbox{ where }
\begin{cases}
y \notin C_1 \cup C_2 \cup C_3 \; \;  &\mbox{if} \; \; (n,3) \neq 1\\
y \neq \{\underline{0}\} \; \; &\mbox{if} \; \;  (n,3)=1.
\end{cases}
\end{equation}
For the structure constants $a, b, c$ of $S$, denote
\[
\kk_0 := \overline{\Qset(a,b,c,y_1,y_2,y_3,y_4)} \subset \kk.
\]
Fix an embedding $\kk_0 \subset \Cset$. Denote by $A^\kk$, $A^{\kk_0}$, and $A^{\Cset}$ the factor algebras
\[
S/( (g-y_4)S + \Sigma_i (z_i - y_i) S)
\]
when the base field is $\kk$, $\kk_0$ and $\Cset$, respectively. Clearly, 
\begin{equation}
\label{ext}
A^\kk \cong A^{\kk_0} \otimes_{\kk_0} \kk, \quad \quad A^\Cset \cong  A^{\kk_0} \otimes_{\kk_0} \Cset.
\end{equation}
\thref{intro2}(1) in the case when the base field is $\Cset$ implies that $A^\Cset \cong M_n(\Cset)$. 
The second isomorphism in \eqref{ext} gives that $A^{\kk_0}$ is a semisimple finite dimensional 
algebra, and thus $A^{\kk_0}$ is a product of matrix algebras over $\kk_0$ (because $\kk_0$ is algebraically closed).
Invoking one more time the second isomorphism in \eqref{ext} gives that $A^{\kk_0} \cong M_n(\kk_0)$. 
The first isomorphism in \eqref{ext} implies that $A^\kk \cong M_n(\kk)$ which completes the 
proof of \eqref{factors} and the proof of \thref{intro2}(1,4)  in the general case.
\qed
%%%%%%%%%%%%%%%%%%

\subsection{\green{Proof of \thref{interm}}} 
\label{proof-last-th}
\green{We begin by providing a discussion of the correspondence between the simple modules over $S$ and the fat point modules over $S$; see, e.g., \cite[Section 3]{Le}. Recall that a {\it fat point module} over $S$ is a 1-critical graded module with multiplicity $>1$. By \cite[Theorem 3.4]{Artin}, we have that all fat point modules over a 3-dimensional Sklyanin algebra $S$  have multiplicity exactly $n/(n,3)$ and thus are $g$-torsionfree. (Indeed, the 1-critical graded modules of $S$ that are $g$-torsion are precisely the point modules of $S$, and these have multiplicity 1.) It is important to point out that since $S$ has Hilbert series $1/(1-t)^3$ we can assume all fat point modules have Hilbert series $d/(1-t)$ with multiplicity $d>1$ up to shift-equivalence. }

\smallskip
\green{On the other hand, let ${\rm Rep}_m(S):={\rm Alg}_\kk(S,{\rm M}_m(\kk))$ be the set of all $m$-dimensional representations over $S$. The algebraic group ${\rm PGL}_m(\kk)\times \kk^\times$ acts on ${\rm Rep}_m(S)$ via $$((T,\lambda).\varphi)(a):=\lambda^iT\varphi(a)T^{-1}$$ for any $\varphi\in {\rm Alg}_\kk(S,{\rm M}_m(\kk))$ and $(T,\lambda)\in {\rm PGL}_m(\kk)\times \kk^\times$, with $a\in S_i$. For simplicity, we write $$\varphi^\lambda:=(1,\lambda).\varphi.$$ It is clear that $\varphi \cong \varphi^\lambda$ if and only if there is some $T\in {\rm PGL}_m(\kk)$ such that $(T,\lambda).\varphi=\varphi$, that is, if $(T,\lambda)$ lies in the stabilizer of $\varphi$.}

\smallskip

\green{To connect the two notions above, a result of Le Bruyn \cite[Proposition~6 and its proof]{Le} and a result of Bocklandt and Symens \cite[Lemma~4]{BS} says that any simple $g$-torsionfree module $V$ over $S$ corresponds to (as simple quotient of) a fat point module $F$ of period $e$ in such a way that
\begin{itemize}
\item $\dim V=de$ with $d={\rm mult}(F)$, and 
\item the stabilizer of $V$ in ${\rm PGL}_{de}(\kk)\times \kk^\times$ is conjugate to the subgroup generated by $(g_\zeta,\zeta)$ with $g_\zeta={\rm diag}(\underbrace{1,\dots,1}_{d},\underbrace{\zeta,\dots,\zeta}_{d},\dots,\underbrace{\zeta^{e-1},\dots,\zeta^{e-1}}_{d})$ and $\zeta$ is a primitive $e$-th root of unity.
\end{itemize}}
\smallskip

\green{Now let us restrict to the case when $n$ is divisible by 3 as in the statement of Theorem~\ref{tinterm}. Let $$\mathfrak{m}:=\mathfrak{m}_p=(z_0-a_0,~z_1-a_1,~z_2-a_2,~g-a_3) \in \maxSpec(Z(S))$$
correspond to a point $p \in (C_1 \cup C_2 \cup C_3) \setminus \{\underline{0}\}$ for some $a_i\in \kk$. Let $V$ be any simple module of $S$ whose central annihilator corresponds to $\mathfrak m$, which can be also considered as a surjective map $\varphi\in {\rm Alg}_\kk(S,{\rm M}_m(\kk))$; here,  $\dim V=m$. From our choice of point $p$, we have $a_3\neq 0$ and $V$ is $g$-torsionfree.}

\smallskip

\green{By the discussion above, $V$ corresponds to some fat point $F$ of period $e$ and multiplicity $n/3$. We claim that $e=1$. By Theorem \ref{tintro2}, $\dim V=e(n/3)<n$ since $\mathfrak{m}$ lies in the singularity locus of $Y$. So $e=1,2$. If $e=2$, then the stabilizer $\varphi$ in ${\rm PGL}_m(\kk)\times \kk^\times$ is conjugate to the subgroup generated by the element $(g_{\zeta},\zeta)$ with $\zeta=-1$. This implies that $V$ is fixed by some $(T,-1)$ in ${\rm PGL}_m(\kk)\times \kk^\times$. Hence $\varphi\cong \varphi^\zeta$ with $\zeta=-1$ and they should share the same central character. But the central character of $\varphi^\zeta$ is given by 
\begin{equation}
\label{CentralChar}
\small{
\begin{array}{rl}
\smallskip
(\varphi^\zeta(z_0),\varphi^\zeta(z_1),\varphi^\zeta(z_3),\varphi^\zeta(g) )&=((-1)^n\varphi(z_0),(-1)^n\varphi(z_1),(-1)^n\varphi(z_3),(-1)^3\varphi(g))\\
\smallskip
& =((-1)^na_0,(-1)^na_1,(-1)^na_3,(-1)^3a_3),
\end{array}
}
\end{equation}
which is not equal to $(a_0,a_1,a_2,a_3)$ when $a_3\neq 0$. Hence, we have that $e=1$. }

\smallskip

\green{As a consequence, $\dim V={\rm mult} (F)=n/3$. Moreover, $\varphi$ has trivial stabilizer, which means $\varphi^\zeta\neq \varphi^{\zeta'}$ whenever $\zeta\neq \zeta'$ again by the above discussion. Therefore, for  $\xi$ is a primitive third root of unity, we have that $\varphi,\varphi^{\xi}$, and $\varphi^{\xi^2}$ are three non-isomorphic irreducible representations whose central characters all correspond to the point $(a_0,a_1,a_2,a_3)$ by \eqref{CentralChar}.}

\smallskip

\green{Finally, denote by $m_1, m_2, \ldots, m_t$ the dimensions of the isomorphism classes of irreducible representations of $S/ (\mm S)$. Accounting for $\varphi$, $\varphi^{\xi}$, $\varphi^{\xi^2}$, we get that $$m_1 = m_2 = m_3 = \dim V=n/3$$ from the discussion above. Now by  \cite[Proposition~4(i)]{Braun}, we obtain that $$n ~=~n/3+n/3+n/3 ~\le~ m_1 + \cdots + m_t ~\leq~ n.$$ (Here we make use of the fact that $Z$ is normal \cite[Corollary on page~2]{Stafford}.) This implies that $t=3$ and every irreducible representation of $S/(\mm S)$ is isomorphic to one of the three representations $\varphi,\varphi^{\xi}$ and $\varphi^{\xi^2}$.
\qed}

\sectionnew{Further Directions and additional results} \label{sec:future}

Application of the techniques above to study the representation theory of other PI elliptic algebras is work in progress. This includes work in preparation for the 4-dimensional Sklyanin algebras (and the corresponding twisted homogeneous coordinate rings) that are module-finite over their center. We believe that our method of employing specializations at levels beyond $N=1$ and then using Poisson geometry will have a wide range of applications to the classifications of the Azumaya loci of the elliptic algebras that are module-finite over their center.

\smallbreak We can combine \thref{intro2} and the recent work of Brown-Yakimov \cite{BrownYakimov} to obtain information for the {\it discriminant ideals}
of the 3-dimensional Sklyanin algebras. 
The role of the discriminant (ideals) first arose in the noncommutative algebra in Reiner's book \cite{Reiner} for the purpose of studying orders and lattices in central simple algebras. There have been several recent applications of noncommutative discriminants including the analysis of automorphism groups of PI algebras (e.g., work of Ceken-Palmieri-Wang-Zhang \cite{CPWZ}) and the Zariski cancellation problem (by Bell-Zhang \cite{BellZhang}). 
A general framework for computing noncommutative discriminants using Poisson geometry was developed in work of Nguyen-Trampel-Yakimov \cite{NTY}.

\smallbreak Towards the goal above, recall that a {\it trace map} on an algebra $A$ is a map $\tr \colon A \to Z(A)$ which is cyclic 
($\tr( xy) = \tr(yx)$ for $x, y \in A$), $Z(A)$-linear, and satisfies $\tr(1) \neq 0$. For a \blue{nonnegative} integer $k$, 
the $k$-th {\em{discriminant ideal}} $D_n(A/Z(A))$ and the $k$-th {\em{modified discriminant ideal}} $MD_k(A/Z(A))$
of $A$ are defined to be the ideals of $Z(A)$ generated by the sets 
\begin{align*}
&\{\det( [\tr(y_i y_j)]_{i,j=1}^k) \mid y_1, \ldots, y_k \in A\} \; \;  \mbox{and}
\\
&\{ \det( [\tr(y_i y'_j)]_{i,j=1}^k) \mid y_1, y'_1, \ldots, y_k, y'_k \in A \},
\; \;  \mbox{respectively}.
\end{align*}

\smallbreak On the other hand, Stafford \cite{Stafford} proved that the PI Sklyanin algebras $S$ are maximal orders in central simple algebras and by \cite[Section~9]{Reiner} 
they admit the so-called {\em{reduced trace maps}}, \blue{which will be denoted by $\tr_{\mathrm{red}}$}.

\blue{As an application of Theorems \ref{tintro2}(4) and \ref{tinterm} we obtain 
a full description of the zero sets of the discriminant ideals of PI 3-dimensional Sklyanin algebras equipped with the reduced trace map.}

\bpr{disc} \blue{For all PI 3-dimensional Sklyanin algebras \blue{$S$} of PI degree $n$,  the zero sets of all discriminant and modified discriminant ideals of $S$ are given by 
\[
\mathbb{V}(D_{k}(S/Z(S)), \tr_{\mathrm{red}})=\mathbb{V}(MD_{k}(S/Z(S)), \tr_{\mathrm{red}})
= 
\begin{cases}
\underline{0}, & k \in [2,n^2] \\
\varnothing, & k =1
\end{cases}
\]
in the case $(n,3) =1$, and by
 \[
\mathbb{V}(D_{k}(S/Z(S)), \tr_{\mathrm{red}})=\mathbb{V}(MD_{k}(S/Z(S)), \tr_{\mathrm{red}})
= 
\begin{cases}
C_1 \cup C_2 \cup C_3, & k \in [n^2/3 +1, n^2] \\
\underline{0}, & k \in [2,n^2/3] \\
\varnothing, & k =1
\end{cases}
\]
in the case $(n,3) \neq1 $. In particular, 
\[
\mathbb{V}(D_{n^2}(S/Z(S)),  \tr_{\mathrm{red}}) = \mathbb{V}(MD_{n^2}(S/Z(S)),  \tr_{\mathrm{red}}) = Y^{\textnormal{sing}}.
\]
}
\epr
\begin{proof} 
\blue{Given $\mm \in \mathrm{maxSpec} Z$, let $\Irr_\mm(S)$ denote the set of isomorphism classes of irreducible representations of $S$ with central annihilator $\mm$. Denote
$$
d(\mm) := \sum_{V \in \Irr_\mm(S)} (\dim_\kk V)^2.
$$
Then \cite[Main Theorem (e)]{BrownYakimov} gives that for all nonnegative integers $k$
\begin{align}
\label{dis}
\mathbb{V}(D_{k}(S/Z(S)), \tr_{\mathrm{red}})&=\mathbb{V}(MD_{k}(S/Z(S)), \tr_{\mathrm{red}})
\\
&=  \Big\{ \mm \in \mathrm{maxSpec} Z \mid
d(\mm)  < k \Big\}.
\nonumber
\end{align}
Theorems \ref{tintro2}(4) and \ref{tinterm} provide a complete classification 
of the irreducible representations of PI 3-dimensional Sklyanin algebras. 
The dimensions of the irreducible representations are also obtained in 
these theorems. Applying these results we get
\begin{equation}
\label{d-funct}
d(\mm_y) = 
\begin{cases}
n^2, & y \in Y \backslash \{ \underline{0} \}
\\
1, & y = \underline{0}
\end{cases}
\quad \mbox{and} \quad
d(\mm_y) = 
\begin{cases}
n^2, & y \in Y \backslash (C_1 \cup C_2 \cup C_3)
\\
n^2/3, & y \in (C_1 \cup C_2 \cup C_3) \backslash \{ \underline{0} \}
\\
1, & y = \underline{0}
\end{cases}
\end{equation}
in the cases  $(n,3) = 1$ and  $(n,3) \neq1$, respectively. Here for $y \in Y$, $\mm_y$ denotes the corresponding maximal ideal of $Z$. 
By combining \eqref{dis} and \eqref{d-funct}, we obtain the statement of the proposition.}
\end{proof}
%\red{MY: I propose to remove this:}
%It is an interesting and important problem to fully describe the discriminant ideals of the elliptic algebras 
%that are module-finite over their center, especially considering the recent applications of noncommutative discriminants discussed above.
%

%%%%%%%%%%%%%%%%%
%%%%%%%%%%%%%%%%%
%%%%%%%%%%%%%%%%%
%%%%%%%%%%%%%%%%%

\sectionnew{Appendix: Examples for $S$ of PI degree 2 and 6} \label{sec:2and6}
In this part, we illustrate \green{Theorems~\ref{tintro}, \ref{tintro2} and \ref{tinterm}} for 3-dimensional Sklyanin algebras of PI degree $n=2$ and $n=6$. We employ the notation of these theorems and of Section~\ref{sec:back} throughout. 

\subsection{PI degree 2} We refer to \cite{ReichWalton} for results on the (representation theory of) 3-dimensional Sklyanin algebra of PI degree 2. Here, $S=S(1,1,c)$ with $c^3 \neq 0, 1, -8$ (cf. \cite[Conjecture~10.37]{AS}), and the center $Z$ is generated by $z_1 = x^2$,  $z_2 = y^2$, $z_3 = z^2$, $g = cy^3 + yxz - xyz -cx^3$, subject to the degree 6 relation:
\[
F = g^2 + \Phi(z_1,z_2,z_3), \quad \text{ with } \Phi = - c^2(z_1^3+z_2^3+z_3^3) - (c^3-4)z_1z_2z_3.
\]
The Poisson $Z$-order structure on $S$, and the induced bracket on $Z$, is described by Theorem~\ref{tintro} with the data above.

\smallbreak Note that $z_1, z_2, z_3$ are second powers of a basis of the generating space $B_1$ of the twisted homogeneous coordinate ring $B=B(E,\mathcal{L},\sigma)=S/gS$, where
\[
\begin{array}{l}
E = \mathbb{V}(\phi) \subseteq \mathbb{P}^2_{[v_1:v_2:v_3]}, \quad \text{ with } \phi = c(v_1^3+v_2^3+v_3^3) - (2+c^3)v_1v_2v_3,\\
\mathcal{L} = \mathcal{O}_{\mathbb{P}^2}(1)|_E,\\ 
\sigma[v_1:v_2:v_3] = [cv_2^2-v_1v_3 : cv_1^2-v_2v_3 : v_3^2 - c^2v_1v_2].
\end{array}
\]
Further, $\{x_1=x,~ x_2=y, ~x_3=z\}$ is a  good basis of $S_1$, and the generators $\{z_1, z_2, z_3\}$ are of the form \eqref{good-u}.

\smallbreak Representation-theoretic results on $S$ are given by Theorem~\ref{tintro2} in the case $(n,3) =1$. Here, $Y = \text{maxSpec}(Z) = \mathbb{V}(F) \subseteq \mathbb{A}^4_{(z_1,z_2,z_3,g)}$,
 which admits an action of the group $\Sigma:= \mathbb{Z}_3 \times \kk^\times$. The singularity locus 
 $Y^{\textnormal{sing}}$ of $Y$ is the origin $\{\underline{0}\}$, and the $\Sigma$-orbits of the symplectic cores of $Y$ are $Y \setminus Y_0,  Y_0 \setminus \{\underline{0}\}$ and $\{\underline{0}\},$ with the first two orbits corresponding to the Azumaya part of $Y$. So, the maximal ideals in $Y \setminus \{\underline{0}\} = [Y \setminus Y_0] ~\cup~ [Y_0 \setminus \{\underline{0}\}]$ are central annihilators of irreducible representations of $S$ of maximum dimension (=2). The maximal ideal corresponding to the origin is the central annihilator of the trivial $S$-module $\kk$. This is consistent with \cite[Theorem~1.3]{Walton}; see \cite[Theorem~7.1]{ReichWalton}.

\subsection{PI degree 6} By \cite[Propositions~1.6 and~5.2]{Walton}, we take the 3-dimensional Sklyanin algebra $S = S(1,-1,-1)$, which has PI degree 6 (cf. \cite[(0.5)]{AS}). A computation shows that $\sigma_{1,-1,-1}^2$ is the permutation $\rho_3$. So by Remark~\ref{runiquegoodbasis}, we have that for $\zeta =e^{2\pi i/3}$, the elements
$$x_1:=x+y+z, \quad x_2:=x+\zeta^2 y+\zeta z, \quad x_3:=x+\zeta y+\zeta^2 z$$ form a good basis of $S$, and the following elements are generators of its center:
\[
z_1 = x_1^6, \quad
z_2 = x_2^6, \quad
z_3 =  x_3^6, \quad
g = x^3-yxz.
\]
Here, $s=2$ and take 
$$u_1 = x_1^2, \quad u_2 = x_2^2, \quad u_3 =  x_3^2,$$ so that $z_i = u_i^3$ for $i =1,2,3$. With the aid of the GBNP package of the computer algebra system GAP \cite{GAP-GBNP}, a calculation shows that the relation $F$ of $Z$ is 
$$F~~= g^6 + 3\ell g^4+3\ell^2g^2+\Phi(z_1,z_2,z_3)$$
where $g^2 + f_3(u_1,u_2,u_3) = 0$ in $Z(S^{(3)})$, and
\[
\hspace{-.05in}\begin{array}{l}
\ell =\frac{1}{108} (z_1+z_2+z_3), \quad \Phi = \ell^3 - \textstyle\frac{1331}{373248} z_1z_2z_3, \quad 
f_3 = \frac{1}{108}\left(u_1^3+u_2^3+u_3^3-\frac{132\zeta^2}{8}u_1u_2u_3\right).
\end{array}
\]

We will now see that representation-theoretic results on $S$ are consistent with  Theorem~\ref{tintro2} in the case $(n,3) \neq 1$. Here, $Y = \text{maxSpec}(Z) = \mathbb{V}(F) \subseteq \mathbb{A}^4_{(z_1,z_2,z_3,g)}$, which admits an action of the group $\Sigma:= \mathbb{Z}_3 \times \kk^\times$. The singularity locus $Y^{\textnormal{sing}}$ of $Y$ is the union of  curves $C_1,  C_2,  C_3$, where 
\[
\begin{array}{rlll}
C_1 &= \{z_1 = -108g^2, ~&z_2 = z_3 = 0, 
~&g \text{ free}\}, \\
C_2 &= \{z_2 = -108g^2, ~&z_1 = z_3 =0, 
~&g \text{ free}\}, \\
C_3 &= \{z_3 = -108 g^2, ~&z_1= z_2 = 0; 
~&g \text{ free}\},
\end{array}
\]
each curve is invariant under dilation \eqref{dilation} (cf. Lemma~\ref{lYsing}). Now for $\gamma \in \kk$, we have that $Y^{\textnormal{sing}}_{\gamma =0}$ is the origin, and that
$$Y^{\textnormal{sing}}_{\gamma \neq 0} = 
\textstyle \left\{\left(-108\gamma^2, ~0, ~ 0\right), ~\left(0, ~-108 \gamma^2,  0\right), ~ \left(0, 0, ~ -108\gamma^2\right)\right\}$$
the union of 3 distinct points.  The $\Sigma$-orbits of the symplectic cores of $Y$ are 
$$Y \setminus (Y_0 \cup C_1 \cup C_2 \cup C_3), \quad \quad (C_1 \cup C_2 \cup C_3) \setminus \{\underline{0}\}, \quad \quad Y_0\setminus  \{\underline{0}\}, \quad \quad  \{\underline{0}\}, $$ 
with the first and third orbits corresponding to the Azumaya part of $Y$. So, the maximal ideals in $[Y \setminus (Y_0 \cup C_1 \cup C_2 \cup C_3)] ~\cup~ [Y_0\setminus  \{\underline{0}\}]$ are central annihilators of irreducible representations of $S$ of maximum dimension (=6). The maximal ideal corresponding to the origin is the central annihilator of the trivial $S$-module $\kk$.

\smallbreak 
%Now we verify Conjecture~\ref{jinterm} for $S(1,-1,-1)$ of PI degree 6. 
\green{Finally, we illustrate Theorem~\ref{tinterm} for the Sklyanin algebra $S(1,-1,-1)$ (of PI degree~6)}. Using the $\Sigma$-action we only need to 
display three non-isomorphic 2-dimensional irreducible representations of $S$ annihilated by $\mathfrak{m} \in \text{maxSpec}(Z)$ corresponding to a point $p$ on $C_1 \setminus \{\underline{0}\}$.
Take $p = \left(-108(4)^2, ~0, ~0, ~4\right) \in C_1$. Then considering the representation $\varphi$ of $S$ (in terms of its good basis), we get that the three representations $\varphi$, $\zeta \varphi$, $\zeta^2 \varphi$ fulfill our goal:
\footnotesize{\[
\begin{array}{lll}
\varphi(x_1) = \begin{pmatrix} 3-i & -1-i \\ 1-i &3+i \end{pmatrix} & 
\varphi(x_2) = \begin{pmatrix} -i & -\zeta-i\zeta^2 \\ \zeta-i\zeta^2 &i \end{pmatrix} & 
\varphi(x_3) = \begin{pmatrix} -i & -\zeta^2-i\zeta   \\ \zeta^2-i\zeta   &i \end{pmatrix}.
\end{array}
\]}

%%%%%%%%%%%
\def\cprime{$'$}

\end{document}